\documentclass{article}

 \usepackage{amsmath}
 \usepackage{amsthm}
 \usepackage{amsfonts}
 \usepackage{amssymb}
 \usepackage{graphicx}

\def\R{\mathbb{R}}

\def\F{\lvert F \rvert}
\def\gnk{G(n,k)}
\def\gn2{G(n,2)}
\def\P{\lvert P \rvert}
\def\npi{\lvert \Pi \rvert}

\def\I{\lvert I \rvert}
\def\tilpi{\tilde{\Pi}}
\def\tili{{\tilde{I}}}
\def\ntilpi{\lvert \tilde{\Pi} \rvert}
\def\ntili{\lvert \tilde{I} \rvert}
\def\tilik{\tilde{I}_k}
\def\ntilik{\lvert \tilde{I}_k \rvert}

\def\nvk{\lvert V_k \rvert}
\def\nvkp{\lvert V_{k,p} \rvert}
\def\nsigpix{\lvert \Sigma_{\pi_0,x} \rvert}
\def\nppix{\lvert P_{\pi_0,x} \rvert}

\theoremstyle{plain}
  \newtheorem{theorem}[subsection]{Theorem}
  \newtheorem{conjecture}[subsection]{Conjecture}

  \newtheorem{proposition}[subsection]{Proposition}
  
  \newtheorem{lemma}[subsection]{Lemma}
  \newtheorem{corollary}[subsection]{Corollary}

\theoremstyle{remark}
  \newtheorem{remark}[subsection]{Remark}
  
  \newtheorem{acknowledgements}[subsection]{Acknowledgements}

\theoremstyle{definition}
  \newtheorem{definition}[subsection]{Definition}

\title{An Incidence Bound for $k$-planes in $F^n$ and a Planar Variant of the Kakeya Maximal Function}

\author{John Bueti\\
\small{University of California, Riverside}\\
\small{Department of Mathematics}\\
\small{jbueti@math.ucr.edu}}
\date{}

\begin{document}
\maketitle

\begin{abstract}
\noindent We discuss a planar variant of the Kakeya maximal
function in the setting of a vector space over a finite field.
Using methods from incidence combinatorics, we demonstrate that
the operator is bounded from $L^p$ to $L^q$ when $1 \leq p \leq
\frac{kn+k+1}{k(k+1)}$ and $1 \leq q \leq (n-k)p'$.
\end{abstract}

\section{Introduction}\label{intro}

The Kakeya conjecture is a long standing open problem in the field
of geometric combinatorics which is concerned with the extent to
which a large direction-separated family of thin tubes can be
compressed into a small space.  There are essentially two
formulations of the conjecture: one geometric and the other
analytic.  In order to state the geometric formulation of the
conjecture, we shall need the following fundamental definition:

\begin{definition}
A set $E \subset \R^n$ is said to be a \emph{Kakeya set} if for
any direction $\xi \in S^{n-1}$, there exists a unit line segment
$l_{\xi}$ parallel to $\xi$ such that $l_{\xi} \subset E$.
\end{definition}

The Kakeya conjecture is concerned with the dimension of such
objects.  Explicitly,

\begin{conjecture}\label{kakeyaconjecture:dim}
If $E \subset \R^n$ is a Kakeya set, then $\dim (E)=n$.
\end{conjecture}

Strictly speaking, Conjecture \ref{kakeyaconjecture:dim} should be
interpreted as three separate conjectures, as one can consider the
Hausdorff dimension (denoted $\dim_H(E)$), lower Minkowski
dimension (denoted $\underline{\dim}(E)$) or upper Minkowski
dimension (denoted $\overline{\dim}(E)$) of such sets; these
different notions of dimension are discussed at great length in
most texts on geometric measure theory, for example
\cite{falconer} and \cite{mattila}.  By definition, one has
\begin{equation*}
\dim_H(E) \leq \underline{\dim}(E) \leq \overline{\dim}(E)
\end{equation*}
for any set $E \subset \R^n$.  Therefore, any progress on the
Hausdorff version of the Kakeya conjecture will immediately imply
progress on the other two versions.  It is not known, however,
whether all three versions of the conjecture are actually
equivalent.  Furthermore, the best known results concerning the
conjecture differ for each of the various notions of dimension
(see \cite{kt:kakeya2} and \cite{katzlabatao} for examples of this
phenomenon).

The analytic formulation of the conjecture concerns the
boundedness of the Kakeya maximal function:

\begin{definition}
Given a function $f \in L^1_{loc}(\R^n)$ and a number $0 < \delta
\ll 1$, we may define the \emph{Kakeya maximal function} as
\begin{align*}
f^*_{\delta} &: S^{n-1} \rightarrow \R\\
f^*_{\delta}(\xi) &=\sup_{a \in \R^n}\frac{1}{\lvert
T_{\delta}(a,\xi)\rvert}\int_{T_{\delta}(a,\xi)}\lvert f(x) \rvert
dx
\end{align*}
where $T_{\delta}(a,\xi)$ denotes the tube of dimensions
$1\times\delta^{n-1}$ centered at the point $a$, oriented in the
direction $\xi$.
\end{definition}

With this definition, we may state the analytic version of the conjecture:

\begin{conjecture}
For all $\epsilon > 0$, and $0<\delta\ll 1$, one has
\begin{equation}
\lVert f^{*}_{\delta}  \rVert_{L^p(S^{n-1},d\sigma)} \leq
C_{\epsilon}\delta^{\frac{n}{p}-1-\epsilon}\lVert f
\rVert_{L^p(\R^n,dx)}
\end{equation}
for $1 \leq p \leq n$, where $d\sigma$ denotes the rotationally
invariant probability measure on the unit sphere.
\end{conjecture}

It is now known, thanks to an observation of Bourgain
\cite{bourgain:kakeya}, that any progress towards the resolution
of the analytic formulation automatically implies progress towards
the geometric formulation.  In particular, if the Kakeya maximal
function is bounded on $L^p$, it follows that Kakeya sets have
Hausdorff (and, hence, upper and lower Minkowski) dimension at
least $p$.

Modern investigations of this conjecture essentially date back to
the 1970's, when Davies \cite{davies} proved that any Kakeya set
in $\R^2$ has Hausdorff dimension (and, hence, upper and lower
Minkowski dimension) $2$.  Later, C\'ordoba \cite{cordoba} proved
the analogous  result for the maximal function in $\R^2$.  Though
these two results essentially resolved all questions concerning
the Kakeya conjecture in $2$ dimensions, the conjecture remains
open in dimensions $3$ and higher.

The $(n,k)$ problem is a variant of the Kakeya problem in which
one replaces lines with $k$-planes.  In order to formally describe
the $(n,k)$ problem, we need a few definitions.

\begin{definition}
A set $E \subset \R^n$ is said to be an $(n,k)$ \emph{set} if for
any $k$-dimensional subspace $\pi \subset \R^n$, there exists a
$k$-dimensional unit cube $Q_{\pi}$ parallel to $\pi$ such that
$Q_{\pi} \subset E$.
\end{definition}

Clearly, an $(n,1)$ set is simply a Kakeya set.  The geometric
conjecture associated with these sets is essentially the same as
the Kakeya conjecture:

\begin{conjecture}
If $E \subset \R^n$ is an $(n,k)$ set, then $\dim (E) = n$.
\end{conjecture}

Once again, dimension can be taken as upper Minkowski, lower
Minkowski or Hausdorff.  One can also define the $(n,k)$ maximal
function:

\begin{definition}
Given positive integers $1 \leq k < n$, let $G(n,k)$ denote the
Grassmannian manifold of all $k$-dimensional subspaces of $\R^n$.
For a function $f \in L^1_{loc}(\R^n)$ and a real number $0 <
\delta \ll 1$, we define the $(n,k)$ \emph{maximal function} as
\begin{align*}
T_{n,k,\delta}f& :G(n,k) \rightarrow \R\\
T_{n,k,\delta}f(\pi)& = \sup_{a \in \R^n}\frac{1}{\lvert
P_{\delta}(a,\pi) \rvert}\int_{P_{\delta}(a,\pi)}\lvert f(x)
\rvert dx
\end{align*}
where $P_{\delta}(a,\pi)$ denotes the $\delta$-neighborhood of the
$k$-dimensional unit cube parallel to $\pi$ and centered at $a \in
\R^n$.
\end{definition}

Using this definition, the analytic version of the $(n,k)$ conjecture is then

\begin{conjecture}
For all $\epsilon > 0$, and $0 < \delta \ll 1$, one has
\begin{equation}
\lVert T_{n,k,\delta}f \rVert_{L^p(G(n,k),d\pi)} \leq
C_{\epsilon}\delta^{\frac{n}{p}-k-\epsilon}\lVert f
\rVert_{L^p(\R^n,dx)}
\end{equation}
for $1 \leq p \leq \frac{n}{k}$, where $d\pi$ denotes the
rotationally invariant probability measure on $G(n,k)$.
\end{conjecture}

Though a large amount of work has been devoted to the $(n,k)$
problems mentioned here, as well as various other $k$-plane
transforms (for example, see the work of Strichartz
\cite{strichartz}, Christ \cite{christ:kplanes} and Drury
\cite{drury:riesz}), the previous results most relevant to this
paper are those due to Oberlin and Stein \cite{o-s}, Themis Mitsis
\cite{mitsis:n2},\cite{mitsis:max}, as well as Wolff's bound for
the Kakeya maximal function \cite{wolff:kakeya}

We investigate the problem in the setting of vector spaces over
finite fields; this slightly unconventional choice of vector
spaces was inspired by a recent paper due to Mockenhaupt and Tao
\cite{mocktao}.  When working in a discrete setting, one is
allowed to focus almost exclusively on the geometric and
combinatorial aspects of the problem.  In particular, problems
which arise in the study of Kakeya-type problems in Euclidean
spaces due to multiplicity of scales (in particular, the
troublesome behavior of $\delta$-tubes which intersect at small
angles) are essentially non-existent in the discrete setting.
Furthermore, working in this setting allows one to import useful
tools from other areas of mathematics, such as analytic number
theory (Gauss sums, Kloosterman sums) and incidence geometry.
There are, of course, a few negative side-effects resulting from
the discretization of the problem. For example, Taylor
approximations of surfaces don't make sense in $F^n$, and there is
no chance to make use of arguments requiring induction on scales.

Before stating the main result, a bit of notation is necessary.
Given a finite field $F$, we let $G(n,k)$ denote the
(Grassmannian) set of all $k$-dimensional subspaces of the vector
space $F^n$.  We will always assume that $F$ is a very large
finite field.  Aside from the fact that the results discussed here
are somewhat trivial when $F$ is small, working with large fields
allows one to develop an intuition regarding the Minkowski
dimension of the geometric objects we will be studying.  When
considering Lebesgue spaces of functions defined on $\gnk$, we
endow this set with a normalized counting measure $d\nu$, so that
for any subset $\Pi \subset G(n,k)$, $\nu(\Pi) =
\F^{-k(n-k)}\lvert \Pi \rvert$ (where $\lvert \cdot \rvert$ simply
denotes cardinality).  The vector space $F^n$ is endowed with a
standard counting measure $dx$.  With this notation, we may define
the $(n,k)$ maximal function operator as
\begin{equation}\label{def:nkmax}
T_{n,k}f(\pi) = \sup_{x \in F^n}\sum_{y \in x + \pi}\lvert f(y)
\rvert
\end{equation}
where $f$ is some real-valued function on $F^n$, and $\pi \in
\gnk$.

Our primary objective here will be to determine a class of
Lebesgue spaces on which this operator is ``bounded''.  Of course,
since we are working in vector spaces over finite fields, we won't
be encountering any divergent integrals.  Our definition of
boundedness will be as follows:

\begin{definition}\label{definition}
Let $T_{n,k}(p \rightarrow q)$ denote the smallest quantity such
that
\begin{equation*}
\lVert T_{n,k}f \rVert_{L^q(d\nu)} \leq T_{n,k}(p \rightarrow
q)\lVert f \rVert_{L^p(dx)}
\end{equation*}
holds for all functions $f: F^n \rightarrow \R$.  We say that
$T_{n,k}$ is \emph{bounded} from $L^p$ to $L^q$ if $T_{n,k}(p
\rightarrow q) \lessapprox 1$ as $\F \rightarrow \infty$.
\end{definition}

The notation $A \lessapprox B$ in the above definition means that
for all $\epsilon > 0$, there exists some constant $C_{\epsilon}$
such that $A \leq C_{\epsilon}\F^{\epsilon}B$.  Similarly, $A
\lesssim B$ indicates $A \leq CB$ for some constant $C$
independent of the field $F$.  Also, we say $A \approx B$
(respectively $A \sim B$) if $A\lessapprox B$ and $B \lessapprox
A$ (respectively $A \lesssim B$ and $B \lesssim A$).

One can immediately observe a few necessary conditions on $p$ and
$q$ by examining certain counterexamples.  For example, observing
the function $f \equiv 1$ shows that one must have $p \leq
\frac{n}{k}$.  Also, considering the function $f = \chi_{\pi}$ for
some $k$-plane $\pi \subset F^n$ will reveal that one must have $q
\leq (n-k)p'$.  It is conjectured that these necessary conditions
are also sufficient, though only partial results are currently
available.

Our main result will be a proof of the following result for the
$(n,k)$ maximal function by means of incidence combinatorial
techniques:

\begin{theorem}\label{mainthm}
Let $n$, $k$ be positive integers such that $2 \leq k \leq n-2$.
The operator $T_{n,k}:L^p(F^n) \rightarrow L^q(\gnk)$ is bounded
when $1 \leq p \leq \frac{kn+k+1}{k(k+1)}$ and $1 \leq q \leq
(n-k)p'$.
\end{theorem}

As an immediate corollary, we have the following result concerning
the geometric version of the $(n,k)$ problem:

\begin{corollary}
Let $E \subset F^n$ be an $(n,k)$ set, where $2 \leq k \leq n-2$
are integers.  Then,
\begin{equation}
\lvert E \rvert \gtrapprox \F^{\frac{kn+k+1}{k+1}}.
\end{equation}
\end{corollary}

\begin{remark}
It should be noted that Theorem \ref{mainthm} is in no way the
best known result for the $(n,k)$ problem.  In the Euclidean case,
Bourgain \cite{bourgain:kakeya} showed that $(n,k)$ sets have
positive $n$-dimensional Lebesgue measure whenever $n \leq2^{k-1}
+2$.  Furthermore, Oberlin (\cite{oberlin:kplanes},
\cite{oberlin:recursive}) has generalized this result to the
corresponding maximal function estimate, and improved the result
for certain values of $n$ and $k$ by making use of recent advances
concerning the X-ray transform (\cite{wolff:x-ray},
\cite{lt:x-ray}, \cite{kt:kakeya2}).  The novelty of Theorem
\ref{mainthm} is in its proof; the proof presented here is
entirely combinatorial, whereas the work of Bourgain and Oberlin
incorporates tools from Fourier analysis, thus leading to stronger
results.
\end{remark}

This result might be interpreted as an extension of the work of
Wolff \cite{wolff:kakeya}, who showed that $(n,1)$ sets have
cardinality $\gtrapprox \F^{\frac{n+2}{2}}$, and Oberlin and Stein
\cite{o-s}, whose work implies that $(n,n-1)$ sets have
cardinality $\sim \F^n$.  A Euclidean version of the case $k=2$
was first proven by Mitsis \cite{mitsis:n2}.  His proof is an
adaptation of the ``hairbrush" construction used by Wolff to
demonstrate a bound for the Kakeya maximal function.  For our
purposes, we shall be following the example set by Mockenhaupt and
Tao \cite{mocktao} who find an alternate proof of Wolff's Kakeya
result by incidence combinatorial methods\footnote{It is, of
course, possible to adapt Wolff's method to provide the same bound
stated in Theorem \ref{mainthm}.  Furthermore, the two methods
are, in a certain sense,``equivalent".  For more details, the
reader is directed to \cite{bueti:diss}.}.

The paper will be organized as follows:  In section \ref{prelim},
we will introduce the main combinatorial tools which will be used
in later computations, and establish a correspondence between
bounds on the number of incidences between points and $k$-planes
and estimates for the $(n,k)$ maximal function as defined in
equation \eqref{def:nkmax}.  Next, in section \ref{reductions}, we
will prove a generalized version of Wolff's ``two-ends reduction".
This reduction will allow us to prove non-trivial incidence bounds
without being hindered by the existence of certain pathological
configurations of points and $k$-planes.  The remaining sections
will be devoted to proving an incidence bound between points and
$k$-planes.  We shall arrive at this bound by estimating the
number of $(k+1)$-simplices arising from a configuration of points
$P$ and $k$-planes $\Pi$ which have vertices from $P$ and faces
from $\Pi$.  The lower bound for the number of such simplices is
computed by an inductive procedure, and is addressed in sections
\ref{firstproof} and \ref{indreg}.  The upper bound for the number
of simplices is computed in section \ref{upperbound}.

\begin{acknowledgements}
I would like to acknowledge Terence Tao for introducing me to this
problem and for his guidance.  Also, I would like to thank Richard
Oberlin for his careful reading of, and helpful comments
concerning this paper.
\end{acknowledgements}


\section{Preliminary Incidence Combinatorial Techniques}\label{prelim}
As the proof of the main theorem in this paper will be largely
combinatorial, our first task will be to establish some machinery
designed to translate incidence combinatorial results into maximal
function estimates, and vice versa.

An important combinatorial tool which will be used extensively the
following set-theoretical version of Cauchy-Schwarz:

\begin{theorem}[Cauchy-Schwarz]\label{stcs}
Let $A$ and $B$ be finite sets with some relation $a \sim b$
between elements $a \in A$ and $b \in B$.  Then,
\begin{equation}
\left\lvert \left\{(a,a',b) \in A \times A \times B:
\begin{array}{l}
a \sim b\\
a' \sim b
\end{array}
\right\} \right\rvert \geq \frac{\lvert \{(a,b) \in A \times B: a
\sim b\} \rvert^2}{\lvert B \rvert}.
\end{equation}
\end{theorem}

\begin{proof}
To see how Theorem \ref{stcs} follows from the traditional
Cauchy-Schwarz inequality, define a function $f$ on the set $B$ as
\begin{align*}
f(b)& = \lvert \{a \in A : a \sim b \}\rvert\\
{}& = \sum_{a \in A}\chi_{a \sim b}.
\end{align*}

Next, using Cauchy-Schwarz, we have
\begin{equation*}
\left(\sum_{B}f(b) \right)^2 \leq \lvert B \rvert \lVert f
\rVert_{L^2(B)}^2,
\end{equation*}
where the $L^2$ norm is computed with respect to counting measure
on the finite set $B$.  Theorem \ref{stcs} then follows by
observing that
\begin{equation*}
\left(\sum_{B}f(b) \right)^2 = \lvert \{(a,b) \in A \times B: a
\sim b\} \rvert^2,
\end{equation*}
and
\begin{equation*}
\lVert f \rVert_{L^2(B)}^2 = \lvert\{(a,a',b) \in A \times A
\times B: a \sim b, a' \sim b \}\rvert.
\end{equation*}
\end{proof}

In practice, the sets $A$ and $B$ will always denote certain
configurations of points, planes and lines, and the relation
$\sim$ will denote some sort of geometric incidence.  For example,
given a configuration of points $P$ and $k$-planes $\Pi$, we may
obtain a lower bound on the number of ``double point-plane
incidences'' in terms of the numbers of point-plane incidences and
planes:
\begin{equation}\label{cauchy-schwarz}
\lvert \{(p,p',\pi) \in P \times P \times\Pi: p,p' \in \pi\}\rvert
\geq \frac{\lvert\{(p,\pi) \in P \times \Pi: p \in
\pi\}\rvert^2}{\npi}.
\end{equation}
Of course, this idea can be generalized into a version of
H\"{o}lder's inequality:
\begin{equation}\label{holder}
\lvert \{(p_1,\dotsc,p_m,\pi) \in P^m \times \Pi: p_i \in \pi\}
\rvert \geq \frac{\lvert \{(p,\pi) \in P \times \Pi: p \in
\pi\}\rvert^m}{\npi^{m-1}}.
\end{equation}

The set appearing on the right hand side of both
\eqref{cauchy-schwarz} and \eqref{holder} is known as the
\emph{incidence set} associated to the points $P$ and $k$-planes
$\Pi$.  As this set will be appearing quite often throughout this
paper, we'll introduce the following notation:
\begin{equation}\label{incidence-def}
I(P,\Pi) := \{(p,\pi) \in P \times \Pi: p \in \pi\}
\end{equation}
Often, the arguments $P$ and $\Pi$ will be suppressed when they
are obvious from the context.

The following correlation between maximal function estimates and
incidence bounds is used in \cite{mocktao} to demonstrate a bound
for the Kakeya maximal function, and we generalize it to suit our
purposes as follows:

\begin{proposition}\label{max-ic}
Given exponents $1 \leq p,q \leq \infty$, the bound $T_{n,k}(p
\rightarrow q) \lessapprox 1$ holds if and only if given any
collection of points $P \subset F^n$ and any direction separated
collection of $k$-planes $\Pi$ contained in $F^n$, the following
incidence bound holds:
\begin{equation}\label{icbound}
\lvert \{(p,\pi) \in P \times \Pi \colon p \in \pi\}\rvert
\lessapprox
\P^{\frac{1}{p}}\npi^{\frac{1}{q'}}\F^{\frac{k(n-k)}{q}}.
\end{equation}
\end{proposition}

Before proceeding with the proof, note that the conjectured best
possible incidence bound (corresponding to the necessary
conditions $p \leq \frac{n}{k}$ and $q \leq (n-k)p'$) is
\begin{equation}\label{icbest}
\lvert I(P,\Pi) \rvert \lessapprox
\P^{\frac{k}{n}}\npi^{\frac{n-1}{n}}\F^{\frac{k(n-k)}{n}}.
\end{equation}
This expression will appear several times throughout the course of
the paper.

\begin{proof}
First, we assume $T_{n,k}(p \rightarrow q) \lessapprox 1$.  Let
$P$ and $\Pi$ be as in the statement of the proposition, and let
$D \subset \gnk$ denote the direction set of $\Pi$ (i.e. each $\pi
\in \Pi$ is a parallel translate of exactly one element of $D$).
Then, the result simply follows from H\"older's inequality:
\begin{align*}
\lvert \{(p,\pi) \in P \times \Pi \colon p \in \pi\}\rvert& \leq \sum_{d \in D}T_{n,k}\chi_{P}(\pi(d))\\
{}& = \F^{k(n-k)}\int_{D}T_{n,k}\chi_{P}d\nu\\
{}& \lessapprox \F^{k(n-k)}
\left(\frac{\npi}{\F^{k(n-k)}}\right)^{\frac{1}{q'}}\P^{\frac{1}{p}}.
\end{align*}

To prove the converse, it suffices (by duality) to show
\begin{equation*}
\sum_{x \in F^n}\int_{\gnk}g(\sigma)\chi_{\sigma +
x_0(\sigma)}(x)f(x)d\nu(\sigma) \lessapprox \lVert g
\rVert_{L^{q'}(d\sigma)} \lVert f \rVert_{L^p(dx)}
\end{equation*}
for any pair of functions $f$ and $g$ (defined on $F^n$ and
$\gnk$, respectively), where $x_0$ is some function which
translates elements of $\gnk$ to affine position. Since out
notation allows us to lose factors of $\log \F$, we may employ the
dyadic pigeonhole principle, and assume that $f = \chi_P$ and $g =
\chi_D$ for sets $P \subset F^n$, and $D \subset \gnk$. After
making these simplifications, we must now show
\begin{equation*}
\frac{1}{\F^{k(n-k)}}\sum_{x \in P}\sum_{\sigma \in
D}\chi_{\sigma+x_0(\sigma)}(x) \lessapprox \left(\frac{\lvert D
\rvert}{\F^{k(n-k)}}\right)^{\frac{1}{q'}}\lvert P
\rvert^{\frac{1}{p}}.
\end{equation*}

Now, simply define the collection of $k$-planes to be $\Pi :=
\{\sigma + x_0:\sigma \in D\}$. The above equation then follows
directly from \eqref{icbound}.
\end{proof}


\section{Avoiding Obstructions To Nontrivial \\
Incidence Bounds}\label{reductions}

In this section, we address an issue which often causes problems
when counting incidences between points and (generally speaking)
algebraic varieties of dimension greater than $1$: it is quite
easy to construct a large set of points $P$ and a large set of
planes $\Pi$ such that every point is contained in every plane.
The classic example of this behavior can be seen in the following
example:

\begin{definition}
A \emph{type-$(2,1)$ degenerate configuration} is a configuration
of a set of points $P$ and a set of 2-planes $\Pi$ such that all
of the points in $P$ lie on some line $l$, and all of the planes
in $\Pi$ contain the line $l$.
\end{definition}

\begin{figure}[htb]
\begin{center}
\includegraphics{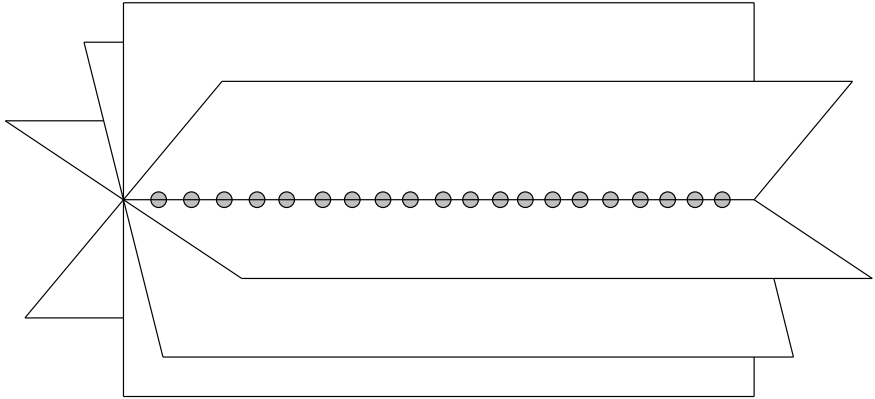}
\caption[A type-$(2,1)$ degenerate configuration.]{A type-$(2,1)$
degenerate configuration.}
\end{center}
\label{21degpic}
\end{figure}

This configuration is notable because the incidence set associated
to it actually attains the worst possible upper bound, $\I = \npi
\P$.  So, it seems that this counterexample should prohibit us
from obtaining any sort of nontrivial upper bound on the size of
the incidence set which would hold for all possible collections of
points and planes.  For higher values of $k$, there are even more
possible degenerate configurations of this type; given a positive
integer $r < k$, one can consider a family of $\F^r$ points which
all lie on some affine $r$-plane $\sigma$, and a direction
separated collection $\F^{(k-r)(n-k)}$ $k$-planes which all
contain $\sigma$.  We shall call such a configuration
\emph{type-$(k,r)$ degenerate}.

Because of the existence of such configurations, it is necessary
to place some (presumably mild) restrictions on distributions of
points and/or planes in question.  In recent work concerning
incidences between points and surfaces, (for example, {\L}aba and
Solymosi \cite{laba-sol}, and Elekes and T\'{o}th \cite{el-toth}),
point sets were assumed to have certain uniformity properties to
prohibit pathological configurations.  For our purposes, however,
we shall work with arbitrary distributions of points in $F^n$, and
exploit the fact that any family of planes in consideration must
be direction separated.

The following proposition acts in the same way as the ``two-ends
reduction'' first used by Wolff \cite{wolff:kakeya} to improve
estimates for the Kakeya maximal function.  We will eliminate the
threat posed by configurations in which points tend to cluster
along low dimensional subsets of $k$-planes from $\Pi$ by
eliminating the possibility of having sparsely populated
$k$-planes.  The actual proof of this new version of the reduction
for the $(n,k)$ problem is quite different from Wolff's; it is
more closely related to Drury's work on the X-ray transform
\cite{drury:xray}.

\begin{proposition}[Generalized Two-ends Reduction]\label{hyperred}
Let $P \subset F^n$ be a collection of points, and $\Pi$ a
direction separated collection $k$-planes in $F^n$ such that
\begin{equation}
\lvert I(P,\Pi) \rvert \lesssim \npi\F^{k-1},
\end{equation}
and
\begin{equation}\label{unifdist}
\lvert P \cap \Pi \rvert \approx \frac{\lvert I(P,\Pi)
\rvert}{\npi} \text{ for each $\pi \in \Pi$}.
\end{equation}

Then, the best possible incidence bound holds:
\begin{equation}\label{2endsconclusion}
\lvert I(P,\Pi) \rvert \lessapprox
\P^{\frac{k}{n}}\npi^{\frac{n-1}{n}}\F^{\frac{k(n-k)}{n}}.
\end{equation}
\end{proposition}

\begin{proof}
The proof of this proposition will be inductive.  We shall first
show that \eqref{2endsconclusion} holds under the assumption $\I
\lesssim \npi$ (hence it is safe to assume $\I \gg \npi$), and
then show that the result holds when $\I \lesssim \npi\F^{r}$ for
any integer $0 \leq r \leq k-1$.

Suppose $\I \lesssim \npi$.  Then, making use of the trivial
estimates $\npi \lesssim \F^{k(n-k)}$ and $\P \geq 1$, we have
\begin{align*}
\I& \lesssim \npi^{\frac{n-1}{n}}\F^{\frac{k(n-k)}{n}}\\
{}& \leq
\P^{\frac{k}{n}}\npi^{\frac{n-1}{n}}\F^{\frac{k(n-k)}{n}}.
\end{align*}
Therefore, we may assume $\I \gg \npi$.

Next, assume that $\npi\F^{r-1} \ll \I \lesssim \npi\F^r$ for some
positive integer $r \leq k-1$.  In order to arrive at the desired
conclusion, we will estimate the size of the following set
\begin{equation}
J_r \equiv J_r(P,\Pi) := \{(p_0,\dotsc,p_r,\pi) \in
P^{r+1}\times\Pi:p_i \in \pi \text{ for each $i$}\}.
\end{equation}
Using H\"older's inequality, we arrive at a lower bound of
\begin{equation}
\lvert J_r \rvert \geq \frac{\I^{r+1}}{\npi^r}.
\end{equation}

In order to compute a corresponding upper bound, we break the set
$J_r$ into a union of disjoint subsets in the following manner:
\begin{equation*}
J_r = \bigcup_{j=0}^r J_r^{(j)}
\end{equation*}
where
\begin{equation*}
J_r^{(j)} := \{(p_0,\dotsc,p_r,\pi) \in J_r: \dim
[p_0,\dotsc,p_r]=j\}
\end{equation*}

Now we estimate each of these subsets separately.  Clearly, we
have
\begin{equation}
\lvert J_r^{(0)} \rvert = \I.
\end{equation}
 Next, when $1 \leq j \leq r-1$, we have $\npi$ choices for the $k$-plane, $\approx \left(\frac{\I}{\npi}\right)^{j+1}$ choices for a $(j+1)$-tuple of points which spans $[p_0,\dotsc,p_r]$, and $\F^{j(r-j)}$ choices for the remaining points.  Therefore, we have
\begin{equation}
\lvert J_r^{(j)} \rvert \lessapprox
\frac{\I^{j+1}}{\npi^j}\F^{j(r-j)}.
\end{equation}

Finally, in order to estimate $\lvert J_{r}^{(r)} \rvert$, we
observe that there are $\sim\P^{r+1}$ choices for the
$(r+1)$-tuple of points, and (because the collection of $k$-planes
is direction separated) there are at most $\sim\F^{(k-r)(n-k)}$
$k$-planes for the collection $\Pi$ which can contain the
$r$-dimensional affine space spanned by the already chosen
$(r+1)$-tuple of points.  Putting together all of these
observations, we have
\begin{equation}\label{2ends:jrest1}
\frac{\I^{r+1}}{\npi^{r}} \leq\lvert J_r \rvert \lessapprox
\P^{r+1}\F^{(k-r)(n-k)} + \I
+\sum_{j=1}^{r-1}\frac{\I^{j+1}}{\npi^j}\F^{j(r-j)}.
\end{equation}
Furthermore, since we are assuming $\I \gg \npi\F^{r-1}$, we have
\begin{equation}
\frac{\I^{r+1}}{\npi^r} \gg \frac{\I^{j+1}}{\npi^j}\F^{j(r-j)}
\end{equation}
whenever $1 \leq j \leq r-1$.  Similarly, $\frac{\I^{r+1}}{\npi^r}
\gg \I$.  Therefore, the first term on the right hand side of
equation \eqref{2ends:jrest1} is dominant, and we have
\begin{equation}
\frac{\I^{r+1}}{\npi^{r}} \lessapprox \P^{r+1}\F^{(k-r)(n-k)}.
\end{equation}
A bit of algebraic manipulation then leads us to the following
incidence bound:
\begin{equation}
\I \lessapprox \P \npi^{\frac{r}{r+1}}\F^{\frac{(k-r)(n-k)}{r+1}}.
\end{equation}

Finally, we do a bit more algebra, and observe that
\begin{align*}
\I& \lessapprox \min\{\P \npi^{\frac{r}{r+1}}\F^{\frac{(k-r)(n-k)}{r+1}}, \npi\F^{r},\P\npi \}\\
{}& \lesssim
\P^{\frac{k}{n}}\npi^{\frac{n-1}{n}}\F^{\frac{k(n-k)}{n}}.
\end{align*}
Since this is true whenever $1 \leq r \leq k-1$, we are done.
\end{proof}

\begin{remark}
It should be noted that in the special case when $\P = \F^{r}$ and
$\npi = \F^{(k-r)(n-k)}$, all four of the quantities $\P\npi$,
$\npi\F^{r}$, $\P\npi^{\frac{r}{r+1}}\F^{\frac{(k-r)(n-k)}{r+1}}$
and $\P^{\frac{k}{n}}\npi^{\frac{n-1}{n}}\F^{\frac{k(n-k)}{n}}$
coincide.  These values for $\P$ and $\npi$ are, of course, the
same values that one would see in the construction of a
type-$(k,r)$ degenerate configuration, as was described earlier in
this section. These configurations, therefore, are examples of
situations in which the worst possible incidence bound coincides
with the best possible incidence bound.
\end{remark}

Now that we have eliminated the potential problems caused by the
existence of sparsely populated $k$-planes, we are in a position
to formulate a more flexible version of Proposition \ref{max-ic}.
In order to make use of Proposition \ref{hyperred} in what
follows, we shall need to establish some notation.  Given a
collection of points $P$ and a direction separated collection of
$k$-planes $\Pi$, it is clear that the following equality holds:
\begin{equation*}
\I = \sum_{\pi \in \Pi}\lvert P \cap \pi \rvert.
\end{equation*}
Therefore, by pigeonholing, there exists a subcollection $\tilpi
\subset \Pi$ such that
\begin{equation}\label{refinedI}
\lvert P \cap \pi \rvert \approx \frac{\I}{\ntilpi} \text{ for
each $\pi \in \tilpi$}.
\end{equation}
The only control we have over the size of the set $\tilpi$ is the
trivial estimate $\ntilpi \leq \npi$.  Furthermore, letting
$\tili$ denote the incidence set $I(P,\tilpi)$, we have
\begin{equation}\label{refinedPi}
\ntili = \sum_{\pi \in \tilpi}\lvert P \cap \pi \rvert \approx \I.
\end{equation}

With this notation in place, we may now state and prove the
following proposition:

\begin{proposition}\label{max-ick}
Let $0 \leq a,b,c \leq 1$ be real numbers such that $(n-k)b +c
\geq 1$.  If the following incidence bound holds
\begin{equation}\label{icboundk}
\ntili \lessapprox \P^a \ntilpi^{1-b}\F^{k(1-c)} +
\P\ntilpi^{\frac{k-1}{k}} + \ntilpi\F^{k-1},
\end{equation}
then the statement $T_{n,k}(p \rightarrow q) \lessapprox 1$ is
true with
$p = \frac{(n-k)b+c}{a}$ and\\
 $q = \min \left\{(n-k)p',\frac{(n-k)b + c}{b}\right\}$.
\end{proposition}

\begin{remark}\label{bados}
Unfortunately, the condition $(n-k)b +c \geq 1$ holds only when $n
\geq k+2$ with the values of $a$, $b$ and $c$ to be computed in
the next section.  This technicality prevents us from using the
methods described here to demonstrate an entirely geometric proof
of the result of Oberlin and Stein.
\end{remark}

\begin{proof}
To start with, if the term $\P \ntilpi^{\frac{k-1}{k}}$ dominates,
then we arrive at the best possible bound:
\begin{equation*}
\ntili \lessapprox
(\P\ntilpi^{\frac{k-1}{k}})^{\frac{k}{n}}(\ntilpi\F^k)^{\frac{n-k}{n}}
=\P^{\frac{k}{n}}\ntilpi^{\frac{n-1}{n}}\F^{\frac{k(n-k)}{n}}.
\end{equation*}
So, we may also assume $\ntili \gg \P\ntilpi^{\frac{k-1}{k}}$.

If the term $\ntilpi\F^{k-1}$ dominates the right-hand side of
equation \eqref{icboundk}, then (since our refinements allow us to
assume equation \eqref{unifdist} holds) Proposition \ref{hyperred}
shows that we obtain the best possible incidence bound. Therefore,
we may assume $\ntili \gg \ntilpi \F^{k-1}$.  To finish the proof,
assume $\ntili \lessapprox \P^a \ntilpi^{1-b}\F^{k(1-c)}$.  Since
we have refined our sets in such a way that $\ntili \approx \I$
and $\ntilpi \leq \npi$, we have $\I \lessapprox \P^a
\npi^{1-b}\F^{k(1-c)}$.  Taking a convex combination (making use
of the assumption $(n-k)b +c \geq 1$) of this estimate with the
trivial estimate $\I \leq \npi\F^k$, and applying Proposition
\ref{max-ic} completes the proof.
\end{proof}


\section{Simplex Construction Part One: The Lower
Bound}\label{firstproof} The object of the next three sections
will be to prove the following incidence bound:

\begin{theorem}
Let $P \subset F^n$ be a collection of points and $\Pi$ a
direction separated collection of $k$-planes contained in $F^n$.
Then,
\begin{equation}\label{ic:main}
\ntili \rvert \lessapprox
\P^{\frac{k(k+1)}{k^2+2k+2}}\ntilpi^{\frac{k^2+k+2}{k^2+2k+2}}\F^{\frac{k(k+1)}{k^2+2k+2}}
+ \P\ntilpi^{\frac{k-1}{k}} + \ntilpi\F^{k-1}
\end{equation}
where the sets $\tili$ and $\tilpi$ are as described in equations
\eqref{refinedI} and \eqref{refinedPi}.
\end{theorem}

Once this incidence bound has been demonstrated, we may apply
Proposition \ref{max-ick} to obtain the desired bound for the
$(n,k)$ maximal function. The proof presented here is inspired by
an argument found in \cite{mocktao} to prove a similar incidence
bound in the case $k = 1$; for their result, they obtain upper and
lower bounds on the number of triangles appearing in the
configuration of points and lines (with sides from the given
collection of lines, and vertices from the given collection of
points).  For our purposes here, we will be obtaining upper and
lower bounds on the number of $(k+1)$-simplices appearing in the
configuration, each with $k+2$ $k$-dimensional faces from the
collection $\tilpi$, and $k+2$ vertices coming from the collection
of points $P$.  The proof of the lower bound will be an induction
on $k$ (making use of the fact that a $(k+1)$-simplex is the cone
of $k$-simplex).  The upper bound will make use of the direction
separatedness of the family of $k$-planes.

We begin by demonstrating a lower bound for the number of
$(k+1)$-simplices arising from collections of points $P$ and
$k$-planes $\tilpi$ satisfying certain hypotheses.  The actual
statement of this lower bound will require quite a bit of
terminology and notation, so we begin with a few definitions.

\begin{definition}
Let $n$ and $k$ be integers such that $1 \leq k \leq n-2$.  Given
a collection of points $P \subset F^n$ and $k$-planes $\tilpi$, we
say that this configuration of points and planes satisfies
hypothesis H1($k$,$P$,$\tilpi$) if
\begin{equation}
\ntili \gg \P\ntilpi^{\frac{k-1}{k}}.
\end{equation}
\end{definition}

\begin{definition}
Let $n$ and $k$ be integers such that $1 \leq k \leq n-2$.  Given
a collection of points $P \subset F^n$ and $k$-planes $\tilpi$, we
say that this configuration of points and planes satisfies
hypothesis H2($k$,$P$,$\tilpi$) if
\begin{equation}
\ntili \gg \ntilpi\F^{k-1}.
\end{equation}
\end{definition}

Observe that these hypotheses are derived from the error terms
from the main incidence bound we aim to demonstrate
\eqref{ic:main}.  As the computation of the lower bound the
collection of simplices will be inductive, we shall need to
investigate the behavior of these hypotheses for varying values of
$k$.  Also, observe that we have not yet mentioned the hypothesis
that the set $\tilpi$ be direction separated; this hypothesis will
not be used until section \ref{upperbound} when computing the
upper bound on the collection of simplices, and its absence from
the lower bound computation greatly simplifies the induction
process.

In order to carry out the induction, we shall need to construct a
few sets from our base sets of $P$, $\tilpi$ and $\tili$.  The
first step towards carrying out these constructions will be
observing that the assumptions $\lvert P \cap \pi \rvert \approx
\frac{\I}{\ntilpi}$ for each $\pi \in \tilpi$ and $\ntili \approx
\I$ permit us to make use of Proposition \ref{hyperred}. This
allows us to assume that $\ntili \gg \ntilpi\F^{k-1}$.
Furthermore, we may use the simple argument found in the proof of
Proposition \ref{max-ick} to assume that $\ntili \gg
\P\ntilpi^{\frac{k-1}{k}}$.  In other words, the statements
H1($k$,$P$,$\tilpi$) and H2($k$,$P$,$\tilpi$) are true.

Next, we construct the following set:
\begin{equation}
{\tilik}' := \{(\pi,p_1,\cdots,p_k) \in \tilpi \times P^k:p_i \in
\pi \text{ for all $i$}, \dim[p_1,\dotsc,p_k]=k-1 \}
\end{equation}
Observe that H\"older's inequality, along with Proposition
\ref{hyperred}, gives us a lower bound on $\lvert {\tilik}'
\rvert$ of
\begin{equation}
\lvert {\tilik}' \rvert \gtrsim \frac{\ntili^k}{\ntilpi^{k-1}}
\end{equation}

This set also needs to be refined a bit; we wish to remove
elements of the set ${\tilik}'$ which are degenerate in the sense
that the $(k-1)$-plane spanned by the $k$-tuple $(p_1,\dotsc,p_k)$
carries a small number of points.  To make this refinement, we
introduce a relation $\sim$ on the set ${\tilik}' \times P$
defined as
\begin{equation}
(\pi,p_1,\dotsc,p_k) \sim q \text{ if $q \in [p_1,\dotsc,p_k]$}.
\end{equation}
Since each $\pi \in \tilpi$ contains $\approx
\frac{\ntili}{\ntilpi}$ points from $P$, one should expect each
$(k-1)$-dimensional slice of $\pi$ to contain roughly
$\frac{\ntili}{\ntilpi\F}$ points from $P$. With this intuition,
we may refine ${\tilik}'$ as follows:
\begin{equation}\label{ik-nondeg}
\tilik := \{i_k \in {\tilik}': \lvert\{q \in P: i_k \sim q \}
\rvert\geq \frac{\ntili}{10\ntilpi\F}  \}
\end{equation}

To see that this refinement is okay, we need only show that the
set ${\tilik}' \smallsetminus \tilik$ is small.  A simple estimate
will take care of this: we have $\ntilpi$ choices for the
$k$-plane, each of these $k$-planes contains $\sim\F^k$
hyperplanes, and each element of the set ${\tilik}' \smallsetminus
\tilik$ can have at most
$\left(\frac{\ntili}{10\ntilpi\F}\right)^k$ $k$-tuples of points
on any of these hyperplanes.  So:
\begin{equation*}
\lvert {\tilik}' \smallsetminus \tilik \rvert \lesssim
\left(\frac{\ntili}{10\ntilpi\F}\right)^k \F^k \ntilpi =
\frac{\ntili^k}{10^k\ntilpi^{k-1}} \ll \lvert {\tilik}' \rvert.
\end{equation*}
This allows us to replace ${\tilik}'$ with $\tilik$ without doing
much harm.

Next, we will take pairs of elements of the set $\tilik$, and
identify them along their $k$-tuples:
\begin{equation}
{V_k}' := \{((\pi_0,p_1,\dotsc,p_k),(\pi,q_1,\dotsc,q_k))\in
\tilik \times \tilik:p_i = q_i \text{ for all $i$}\}.
\end{equation}
An application of Cauchy-Schwarz gives us a lower bound on the
size of this set:
\begin{equation}
\lvert {V_k}' \rvert \geq \frac{\ntilik^2}{\P^k}
\end{equation}

This set also needs to be refined; we wish to remove elements from
${V_k}'$ which are degenerate in the sense that $\pi_0 = \pi$.
Such an element of ${V_k}'$ is, in fact, and element of $\tilik$,
so this refinement will be okay if we can show that $\lvert {V_k}'
\rvert \gg \ntilik$.  In order to show this, we make use of
hypothesis H1($k$,$P$,$\tilpi$):
\begin{equation*}
\lvert {V_k}' \rvert \geq \frac{\ntilik^2}{\P^k} \gtrsim \ntilik
\frac{\ntili^k}{\P^k\ntilpi^{k-1}} \gg \ntilik.
\end{equation*}
So, if we define the set $V_k$ as
\begin{equation}\label{def:vk}
V_k := \{(\pi_0,\pi,p_1,\dotsc,p_k) \in {V_k}':\pi_0 \neq \pi\},
\end{equation}
then we have $\nvk \gtrsim \lvert {V_k}' \rvert \gtrsim
\frac{\ntilik^2}{\P^k}$.

Next, to each element of $V_k$ we wish to add a point from $P$
which lives in $\pi \smallsetminus [p_1,\dotsc,p_k]$.
\begin{equation}\label{def:vkp}
V_{k,p} := \{((\pi_0,\pi,p_1,\dotsc,p_k),x) \in V_k \times P: x
\in P \cap (\pi \smallsetminus [p_1,\dotsc,p_k])\}
\end{equation}
Since there are $\approx \frac{\ntili}{\ntilpi}$ points on each
plane, and they do not cluster along low dimensional subspaces by
Proposition \ref{hyperred}, we have the following lower bound
\begin{equation}
\nvkp \gtrapprox \nvk\frac{\ntili}{\ntilpi}.
\end{equation}

One last refinement is needed before we can proceed with the
construction of simplices.  For a given pair $(\pi_0,x) \in \tilpi
\times P$ such that $x \notin \pi_0$, we define a function
$f(\pi_0,x)$ as
\begin{equation*}
f(\pi_0,x) := \lvert\{(\pi,p_1,\dotsc,p_k) \in
\tilik:((\pi_0,\pi,p_1,\dotsc,p_k),x) \in V_{k,p}\}\rvert.
\end{equation*}
Then,
\begin{equation*}
\sum f(\pi_0,x) = \nvkp \gtrapprox \nvk\frac{\ntili}{\ntilpi},
\end{equation*}
so we may pigeonhole this sum to find a family of plane-point
pairs $(\pi_0,x)$ \begin{equation}\label{def:d} D(P,\tilpi) :=
\{(\pi_0,x) \in \tilpi \times P: x \notin \pi_0, f(\pi_0,x)
\gtrapprox \nvk\frac{\ntili}{\ntilpi \lvert D(P,\tilpi) \rvert}\}.
\end{equation}
As was the case for the parameter $\ntilpi$, we have no control
over the parameter $\lvert D(P,\tilpi) \rvert$ other than the
trivial upper bound $\lvert D(P,\tilpi) \rvert\leq \P \ntilpi$.

\begin{remark}
It is somewhat disconcerting that the value of the function
$f(\pi_0,x)$ may be very small on our set $D(P,\tilpi)$.  This
potential problem, however, will be ruled out in Section
\ref{indreg}.
\end{remark}

In summary, we have made the following constructions:
\begin{align*}
\tilik& = \left\{(\pi,p_1,\dotsc,p_k) \in \tilpi \times P^k:
\begin{array}{l}
p_i \in \pi \text{ for each $i$},\\
\dim[p_1,\dotsc,p_k]=k-1\\
\lvert P \cap [p_1,\dotsc,p_k] \rvert \geq
\frac{\ntili}{10\ntilpi\F}
\end{array} \right\}\\
{}& \quad \ntilik \gtrsim \frac{\ntili^{k}}{\ntilpi^{k-1}}\\
V_k& = \left\{((\pi_0,p_1,\dotsc,p_k),(\pi,q_1,\dotsc,q_k))\in
\tilik \times \tilik:
\begin{array}{l}
p_i = q_i \text{ for each $i$}\\
\pi \neq \pi_0
\end{array}
\right\}\\
{}& \quad \lvert V_k \rvert \gtrsim \frac{\ntilik^2}{\P^k}\\
V_{k,p}& = \{((\pi_0,\pi,p_1,\dotsc,p_k),x) \in V_k \times P: x \in P \cap (\pi \smallsetminus [p_1,\dotsc,p_k])\}\\
{}& \quad \lvert V_{k,p} \rvert \gtrapprox \vert V_k \rvert \frac{\ntili}{\ntilpi}\\
f(\pi_0,x)& = \lvert\{(\pi,p_1,\dotsc,p_k) \in \tilik:((\pi_0,\pi,p_1,\dotsc,p_k),x) \in V_{k,p}\}\rvert\\
{}& \quad \sum f(\pi_0,x) = \nvkp \gtrapprox \nvk\frac{\ntili}{\ntilpi}\\
D(P,\tilpi)& = \{(\pi_0,x) \in \tilpi \times P: x \notin \pi_0, f(\pi_0,x) \gtrapprox \nvk\frac{\ntili}{\ntilpi \lvert D(P,\tilpi) \rvert}\}\\
{}& \quad \lvert D(P,\tilpi) \rvert \leq \P\ntilpi
\end{align*}

With this notation established, we may finally begin to compute a
lower bound on the number of simplices.

\begin{lemma}\label{lowerlemma}
Given an arrangement of points $P$ and $k$-planes $\tilpi$
satisfying hypotheses H1($k$,$P$,$\tilpi$) and
H2($k$,$P$,$\tilpi$), let $S_k(P,\tilpi)$ denote the set of
$(k+1)$-simplices in $F^n$ with faces from $\tilpi$ and vertices
from $P$. Then,
\begin{equation}\label{indlemma}
\lvert S_k(P,\tilpi) \rvert \gtrapprox \nvk^{k+1}
\frac{\ntilpi^{k^2-2k-2}}{\P^k\ntili^{k^2-k-2}}.
\end{equation}
where the set $V_k$ is as defined in equation \eqref{def:vk}.
\end{lemma}

\begin{remark}\label{logloss}
Heuristically, one can easily compute this lower bound
on\linebreak
 $\lvert S_k(P,\Pi) \rvert$.  Since a point from $P$ and a $k$-plane from $\Pi$ are incident with probability $\frac{\I}{\P\npi}$, and a simplex consists of $(k+1)(k+2)$ incidences amongst $(k+2)$ points and $(k+2)$ $k$-planes, it follows that
\begin{equation*}
\lvert S_k(P,\Pi) \rvert \gtrsim \P^{k+2}\npi^{k+2}\left(\frac{\I}{\P\npi}\right)^{(k+1)(k+2)}
\end{equation*}
under the unrealistic assumption that all events of point-plane
incidence are independent.  To validate this heuristic rigorously,
we shall make several refinements in order to ensure a certain
amount of geometric regularity regarding the distribution of
points within each $k$-plane, hence creating a logarithmic loss in
the end result.  If a new procedure for demonstrating the
heuristically obvious lower bound can be established, then the
incidence bound might be slightly improved (i.e. we can replace
the symbol $\lessapprox$ with $\lesssim$).
\end{remark}

The proof of this lemma will be a somewhat complicated induction
on $k$.  For the sake of clarity, it will be useful to briefly
outline the proof before proceeding.

\begin{definition}
Given a collection of points $P$ and $k$-planes $\tilpi$, let
C1($k$,$P$,$\tilpi$) denote the statement that the conclusion of
lemma \ref{lowerlemma} holds.
\end{definition}

\begin{definition}
Given a collection of points $P$ and $k$-planes $\tilpi$, let
C2($k$,$P$,$\tilpi$) denote the statement that
\begin{equation}
\lvert S_k(P,\tilpi) \rvert \gtrapprox
\frac{\I^{(k+1)(k+2)}}{\P^{k(k+2)}\ntilpi^{k(k+2)}}.
\end{equation}
\end{definition}

With this terminology established, we may describe the inductive
procedure.  First, we will prove that C2($0$,$R$,$\Sigma$) is true
for any collection of points $R$ and ``affine $0$-planes'' (i.e.
points) in an ambient space of any dimension.  Though the case
$k=0$ is a bit absurd, it will work (formally) as the base of the
induction.  The reader who is unsatisfied with this may choose to
refer to the arguments found in \cite{mocktao}, and begin the
induction at $k=1$. Next, we make the following essentially
trivial observation:

\begin{lemma}
Let $k$ be an integer, $k \geq 1$.  Given a collection of points
$R$, and a collection of $k$-planes $\Sigma$ (living in an ambient
space of any dimension strictly larger than $k$), we have the
following implication:
\begin{equation}
\{\text{H1($k$,$R$,$\Sigma$) $\&$ H2($k$,$R$,$\Sigma$)}\}
\Rightarrow \{\text{C1($k$,$R$,$\Sigma$)} \Rightarrow
\text{C2($k$,$R$,$\Sigma$)}\}.
\end{equation}
\end{lemma}

\begin{proof}
The proof of this sub-lemma is simply a computation making use of
bounds we have already demonstrated.  Use hypotheses
H1($k$,$R$,$\Sigma$) and H2($k$,$R$,$\Sigma$) to assume all of the
regularity conditions needed for the estimates
\begin{align}
\lvert V_k(R,\Sigma) \rvert& \gtrsim \frac{\ntilik^2}{R^k}\\
\lvert I_k(R,\Sigma) \rvert& \gtrsim \frac{\lvert
I(R,\Sigma)\rvert^k}{\lvert \Sigma \rvert^{k-1}}.
\end{align}
where the sets $V_k(R,\Sigma)$ and $I_k(R,\Sigma)$ are defined as
in equations \eqref{def:vk} and \eqref{ik-nondeg} (respectively).

Then, simply insert these bounds into equation \eqref{indlemma}.
\end{proof}

The final ingredient needed to prove lemma \ref{lowerlemma} is a
statement the C2($k-1$,$P_{\pi_0,x}$,$\Sigma_{\pi_0,x}$)
$\Rightarrow$ C1($k$,$P$,$\tilpi$) (under suitable hypotheses),
where the sets $P_{\pi_0,x}$ and $\Sigma_{\pi_0,x}$ are
collections of points and $(k-1)$-planes (respectively) to be
defined later.  With this argument in place, the proof will then,
schematically, look like the following:
\begin{equation*}
\text{C2($0$)} \Rightarrow \text{C1($1$)} \Rightarrow
\text{C2($1$)} \Rightarrow \text{C1($2$)} \Rightarrow
\text{C2($2$)} \dotsm
\end{equation*}

Unfortunately, there are two technicalities we must dispense with
before beginning with the above program.  The method described
will not work unless we can ensure that simplices formed at every
level of the induction are non-degenerate.  In order to avoid
potential non-degeneracies, we must show that the hypotheses
H1($k$,$P$,$\tilpi$) and H2($k$,$P$,$\tilpi$) imply a new pair of
hypotheses H1($k-1$,$Q$,$\Sigma$) and H2($k-1$,$Q$,$\Sigma$) for
suitably defined sets of points $Q$, and $(k-1)$-planes $\Sigma$
(and so on for $(k-2)$-planes, $(k-3)$-planes, etc...).  In order
to simplify the exposition, we will assume, for the moment, that
the regularity hypotheses hold at every level of the induction.
Then, in Section \ref{indreg}, we will verify that this is indeed
the case.

\begin{proof}[Proof of Lemma \ref{lowerlemma}] To begin the proof,
we let $R$ denote a collection of points, and $\Sigma$ a
collection of affine $0$-planes in some ambient space (of any
finite dimension $\geq 1$).  We wish to show that the conclusion
C2($0$,$R$,$\Sigma$) holds; this essentially means that we must
show
\begin{equation}
\lvert S_0(R,\Sigma) \rvert \gtrapprox \lvert I(R,\Sigma)
\rvert^2.
\end{equation}

This estimate, however, is trivial after observing that
$I(R,\Sigma) = R \cap \Sigma$, and $S_0(R,\Sigma)$ denotes the set
of all line segments with distinct endpoints in the set
$I(R,\Sigma)$.

In order to begin the inductive part of the proof, we will need to
introduce a family of $(k-1)$-planes.  Given a set of points $P$,
and $k$-planes $\tilpi$ satisfying H1($k$,$P$,$\tilpi$) and
H2($k$,$P$,$\tilpi$), fix an element $(\pi_0,x) \in D$ (as defined
in equation \eqref{def:d}).  For this particular pair, we define a
family of $(k-1)$-planes $\Sigma_{\pi_0,x}$ defined as:
\begin{equation}\label{def:sigmasub}
\Sigma_{\pi_0,x} := \left\{\sigma \subset \pi_0:
\begin{array}{l}
\exists (\pi,p_1,\dotsc,p_k) \in \tilik \text{ for which}\\
\sigma = [p_1,\dotsc,p_k], ((\pi_0,\pi,p_1,\dotsc,p_k),x) \in
V_{k,p}
\end{array}\right\}
\end{equation}

Also, we need to identify a class of points
\begin{equation}\label{def:psub}
P_{\pi_0,x} := P \cap \pi_0
\end{equation}
which will interact with the set $\Sigma_{\pi_0,x}$.  In Section
\ref{indreg} we will verify that there exists a set
$\tilde{D}(P,\tilpi) \subset D(P,\tilpi)$ such that $\lvert
\tilde{D}(P,\tilpi) \rvert \geq \frac{1}{2}\lvert
D(P,\tilpi)\rvert$ and the hypotheses
H1($k-1$,$P_{\pi_0,x}$,$\Sigma_{\pi_0,x}$) and
H2($k-1$,$P_{\pi_0,x}$,$\Sigma_{\pi_0,x}$) are true for each pair
$(\pi_0,x) \in \tilde{D}(P,\tilpi)$.  Therefore, we may assume
that the configurations of points and planes
$(P_{\pi_0,x},\Sigma_{\pi_0,x})$ induced from this large subclass
enjoy the same regularity properties as $(P,\tilpi)$.

Next, we observe that the set of $(k+1)$-simplices in $S_k$ which
contain the $k$-plane $\pi_0$ as a face and the point $x$ as a
vertex are in one to one correspondence with the set of
$k$-simplices contained in $\pi_0$ whose faces belong to
$\Sigma_{\pi_0,x}$ and vertices belong to $P_{\pi_0,x}$.  In other
words, each such $(k+1)$-simplex is seen as the cone of some
$k$-simplex contained in $\pi_0$.  So, once we have a lower bound
on the number of $k$-simplices contained in $\pi_0$, we will
obtain a lower bound on the number of $(k+1)$-simplices containing
$\pi_0$ and $x$.

We formalize this as follows:  Let $J =
I(P_{\pi_0,x},\Sigma_{\pi_0,x})$ denote the incidence set arising
from the collection of points $P_{\pi_0}$ and the collection of
$(k-1)$-planes $\Sigma_{\pi_0,x}$.  Also, define the set
\begin{equation*}
J_{k,\pi_0,x} := \left\{(\sigma,p_1,\dotsc,p_k) \in
\Sigma_{\pi_0,x} \times P_{\pi_0,x}^k:
\begin{array}{l}
p_i \in \sigma \text{ for all $i$},\\
\dim[p_1,\dotsc,p_k]=k-1
\end{array} \right\}
\end{equation*}

Next, we observe that (by the construction of the sets
$\Sigma_{\pi_0,x}$ and $P_{\pi_0,x}$) there exists a one-to-one
correspondence between the set $J_{k,\pi_0,x}$ and the set
$\{(\pi,p_1,\dotsc,p_k) \in \tilik:((\pi_0,\pi,p_1,\dotsc,p_k),x)
\in V_{k,p}\}$.  So, we have the following:
\begin{align*}
\frac{\lvert J \rvert^k}{\nsigpix^{k-1}}& \sim \lvert J_{k,\pi_0,x} \rvert\\
{}& = \lvert \{(\pi,p_1,\dotsc,p_k) \in \tilik:((\pi_0,\pi,p_1,\dotsc,p_k),x) \in V_{k,p}\} \rvert\\
{}& \gtrapprox \nvk \frac{\ntili}{\ntilpi\lvert D(P,\tilpi)
\rvert}.
\end{align*}

Now we invoke the inductive hypothesis: letting $S_{k-1,\pi_0,x}$
denote the set of $k$-simplices contained in $\pi_0$ with faces
from $\Sigma_{\pi_0,x}$ and vertices from $P_{\pi_0,x}$, we assume
the statement C2($k-1$,$P_{\pi_0,x}$,$\Sigma_{\pi_0,x}$) is true.
Hence,
\begin{equation}
\lvert S_{k-1,\pi_0,x} \rvert \gtrapprox \frac{\lvert J
\rvert^{k(k+1)}}{\nsigpix^{(k-1)(k+1)}\nppix^{(k-1)(k+1)}}.
\end{equation}
Next, we use the previous computation to reinterpret this lower
bound as:
\begin{equation}
\lvert S_{k-1,\pi_0,x} \rvert \gtrapprox \left(\nvk
\frac{\ntili}{\ntilpi \lvert D(P,\tilpi)
\rvert}\right)^{k+1}\frac{1}{\nppix^{(k-1)(k+1)}}
\end{equation}

Now let $S_{k,\pi_0,x}$ denote the set of $(k+1)$-simplices which
contain the $k$-plane $\pi_0$ as a face and the point $x$ as a
vertex.  Recall that every simplex in $S_{k,\pi_0,x}$ is the cone
of a simplex from $S_{k-1,\pi_0,x}$.  So,
\begin{equation}
\lvert S_{k,\pi_0,x} \rvert \gtrapprox \left(\nvk
\frac{\ntili}{\ntilpi \lvert D(P,\tilpi)
\rvert}\right)^{k+1}\frac{1}{\nppix^{(k-1)(k+1)}}.
\end{equation}
Next observe that for any $\pi_0$, $\nppix \approx
\frac{\ntili}{\ntilpi}$.  So, we have
\begin{equation}
\lvert S_{k,\pi_0,x} \rvert \gtrapprox \left(\nvk
\frac{\ntili}{\ntilpi \lvert D(P,\tilpi)
\rvert}\right)^{k+1}\left(\frac{\ntilpi}{\ntili}\right)^{(k-1)(k+1)}.
\end{equation}

To finish the computation, simply observe that this bound holds
for any pair $(\pi_0,x) \in D(P,\tilpi)$.  So, simply multiply the
previous bound by $\lvert D(P,\tilpi) \rvert$, and then insert the
trivial bound $\lvert D(P,\tilpi) \rvert \leq \P \ntilpi$.  This
finishes the proof of the lemma.
\end{proof}


\section{Simplex Construction Part Two: Regularity Hypotheses}\label{indreg}
The material discussed in this section should essentially be
viewed as an appendix to the previous section.  In particular, we
wish to verify that the simplices constructed in the previous
section are not degenerate in the sense that their faces (edges,
vertices, hyper-edges, etc.) do not collapse onto each other.
This nondegeneracy can be demonstrated by showing that the
hypotheses H1($r$,$Q$,$\Sigma$) and H2($r$,$Q$,$\Sigma$) hold for
most of the arrangements of points and $r$-planes $(Q,\Sigma)$
appearing throughout the induction.

Before beginning the process of verifying these hypotheses, it
will be useful to formally describe the construction of the
collections of points and $r$-planes we shall be working with.
The procedure is virtually identical to the one seen in the
previous section for constructing the collections $P_{\pi_0,x}$
and $\Sigma_{\pi_0,x}$. Let $Q$ be a collection of points and
$\Sigma$ a collection of $r$-planes embedded in some
$(r+1)$-dimensional space, and assume that the hypotheses
H1($r$,$Q$,$\Sigma$) and H2($r$,$Q$,$\Sigma$) have been verified.
Working with these sets, we may (as we have already done with
$k$-planes) construct the set $V_{r,p}(Q,\Sigma)$:
\begin{equation}
V_{r,p}(Q,\Sigma) := \left\{
\begin{array}{l}
(\sigma_0,\sigma,p_1,\dotsc,p_r,x)\\
\in \Sigma^2 \times Q^{r+1}
\end{array}:
\begin{array}{l}
[p_1,\dotsc,p_r] \subset \sigma_0 \cap \sigma\\
\dim[p_1,\dotsc,p_r] = r-1\\
x \in \sigma \smallsetminus \sigma_0, \sigma \neq \sigma_0\\
\lvert[p_1,\dotsc,p_r] \cap Q\rvert \geq \frac{\lvert I(Q,\Sigma)
\rvert}{10\lvert \Sigma \rvert\F}
\end{array}
\right\}.
\end{equation}
Since we're assuming the hypotheses H1($r$,$Q$,$\Sigma$) and
H2($r$,$Q$,$\Sigma$) hold, we may compute a lower bound for the
size of this set as
\begin{equation}
\lvert V_{r,p} \rvert \gtrsim \frac{\lvert
I(Q,\Sigma)\rvert^{2r+1}}{\lvert Q \rvert^r
\lvert\Sigma\rvert^{2r-1}}.
\end{equation}

Next, we define a function $f_r$ on the collection of all pairs
$(\sigma_0,x) \in \Sigma \times Q$ such that $x \notin \sigma_0$
as
\begin{equation}
f_r(\sigma_0,x) := \lvert\{(\sigma,p_1,\dotsc,p_r) \in \Sigma
\times Q^r: (\sigma_0,\sigma,p_1,\dotsc,p_r,x)\in
V_{r,p}(Q,\Sigma)\}\rvert.
\end{equation}
Clearly, we have
\begin{equation*}
\sum_{(\sigma_0,x)}f_r(\sigma_0,x) \sim \lvert V_{r,p}(Q,\Sigma)
\rvert.
\end{equation*}
So, we may once again use a dyadic pigeonholing argument to find a
nice collection $D(Q,\Sigma) \subset Q \times \Sigma$ such that
\begin{equation}\label{frestimate}
f_r(\sigma_0,x) \approx \frac{\lvert V_{r,p}(Q,\Sigma)
\rvert}{\lvert D(Q,\Sigma)\rvert} \text{ for all $(\sigma_0,x) \in
D(Q,\Sigma)$}.
\end{equation}

Now, we are in a position to define the relevant collections of
$(r-1)$-planes necessary for the construction.  Given a pair
$(\sigma_0,x) \in D(Q,\Sigma)$, We define a collection of
$(r-1)$-planes embedded in $\sigma_0$ as
\begin{equation}
\Gamma_{\sigma_0,x} := \left\{\gamma \subset \sigma_0:
\begin{array}{l}
\exists (\sigma,p_1,\dotsc,p_r) \in \sigma \text{ for which
$\gamma = [p_1,\dotsc,p_r]$},\\
((\sigma_0,\sigma,p_1,\dotsc,p_r),x) \in V_{r,p}(Q,\Sigma)
\end{array}\right\}.
\end{equation}
Also, we may simply define the associated collection of points as
\begin{equation}
Q_{\sigma_0,x}=Q \cap \sigma_0.
\end{equation}

As it turns out, the verification of H2 at each level of the
induction is far more elementary than the verification of H1; in
fact, H2 essentially holds by definition.  We state this formally
as follows:

\begin{lemma}
Let $(P,\tilpi)$ be an arrangement of points and $k$-planes
embedded in some ambient space as described in equations
\eqref{refinedPi} and \eqref{refinedI}.  Given a pair $(\pi_0,x)
\in D(P,\tilpi)$, define collections of points $P_{\pi_0,x}$ and
$(k-1)$-planes $\Sigma_{\pi_0,x}$ embedded in $\pi_0$ as defined
in equations \eqref{def:psub} and \eqref{def:sigmasub}
respectively.  Then, we have the following:
\begin{equation}
\text{H2}(k,P,\tilpi) \Rightarrow
\text{H2}(k-1,P_{\pi_0,x},\Sigma_{\pi_0,x})
\end{equation}
for each of the pairs $(x,\pi_0)$.  Furthermore, this regularity
holds at all lower levels of the induction in the sense that
(using the notation already established in this section)
\begin{equation}
\text{H2}(r,Q,\Sigma) \Rightarrow
\text{H2}(r-1,Q_{\pi_0,x},\Gamma_{\sigma_0,x})
\end{equation}
for each of the pairs $(\sigma_0,x)$.
\end{lemma}

\begin{proof}
First, let us recall the definitions of the sets $P_{\pi_0,x}$ and
$\Sigma_{\pi_0,x}$:
\begin{align*}
P_{\pi_0,x}& = P \cap \pi_0,\\
\Sigma_{\pi_0,x}& = \left\{\sigma \subset \pi_0:
\begin{array}{l}
\exists (\pi,p_1,\dotsc,p_k) \in \tilik \text{ for which}\\
\sigma = [p_1,\dotsc,p_k], ((\pi_0,\pi,p_1,\dotsc,p_k),x) \in
V_{k,p}
\end{array}\right\}.
\end{align*}
Next, let us recall the definition of $\tilik$:
\begin{equation*}
\tilik = \left\{(\pi,p_1,\dotsc,p_k) \in \tilpi \times P^k:
\begin{array}{l}
\dim[p_1,\dotsc,p_k] = k-1,\\
\lvert [p_1,\dotsc,p_k] \cap P \rvert \geq
\frac{\ntili}{10\ntilpi\F}
\end{array}
\right\}.
\end{equation*}

Suppose that H2($k-1$,$P_{\pi_0,x}$,$\Sigma_{\pi_0,x}$) failed for
some pair $(\pi_0,x) \in D$; this would imply
\begin{equation}
\lvert I(P_{\pi_0,x},\Sigma_{\pi_0,x})\rvert \lesssim \lvert
\Sigma_{\pi_0,x} \rvert \F^{k-2}.
\end{equation}

This, however, is a contradiction to the assumption that
H2($k$,$P$,$\Pi$) holds, as the above definitions imply the
following string of inequalities:
\begin{equation*}
\frac{\ntili}{10\ntilpi\F} \lesssim \frac{\lvert
I(P_{\pi_0,x},\Sigma_{\pi_0,x})\rvert}{\lvert \Sigma_{\pi_0,x}
\rvert} \lesssim \F^{k-2}
\end{equation*}
or,
\begin{equation}
\ntili \lesssim \ntilpi \F^{k-1}.
\end{equation}

We omit the proof for the remaining levels of the induction, as it
is virtually identical to what we have already shown.
\end{proof}

In order to show that property H1 is inherited at all stages of
the induction, one needs to do a bit more work.  The reason for
this is that property H1 does not simply pass freely from level
$k$ to level $k-1$ by definition (as was essentially the case for
property H2).  Loosely speaking, the statement
H1($k$,$P$,$\tilpi$) states that most maximal-rank $k$-tuples of
points in the set $P$ are incident to many $k$-planes from the set
$\tilpi$ (in this case $k=1$, this general property is usually
referred to as ``bilinearity'').  To illustrate the difficulty,
consider the case $k = 2$, $n \geq 4$.  Choose a pair $(\pi_0,x)
\in D$, and construct the sets $P_{\pi_0,x}$ and
$\Sigma_{\pi_0,x}$ in the usual fashion.  We would like to say
that the hypothesis H1($1$,$P_{\pi_0,x}$,$\Sigma_{\pi_0,x}$) is
satisfied, or (more informally) that most of the points in
$P_{\pi_0,x}$ are incident to many lines from $\Sigma_{\pi_0,x}$.
A natural attempt to prove such a statement is to observe that, by
property H1($2$,$P$,$\tilpi$), most pairs of points in $P$ are
incident to many planes from $\Pi$.  In particular, given a point
$p \in P_{\pi_0,x} \subset P$, there are many planes incident to
the pair $(p,x) \in P^2$.  Unfortunately, since $n \geq 4$, these
$2$-planes are not obligated to intersect $\pi_0$ in a line, hence
not contributing to the set $\Sigma_{\pi_0,x}$.  Therefore, one
must come up with a different means by which to verify
H1($1$,$P_{\pi_0,x}$,$\Sigma_{\pi_0,x}$).

\begin{lemma}
Let $(P,\tilpi)$ be an arrangement of points and $k$-planes
embedded in some ambient space as described in equations
\eqref{refinedPi} and \eqref{refinedI}.  Given a pair $(\pi_0,x)
\in D(P,\tilpi)$, define collections of points $P_{\pi_0,x}$ and
$(k-1)$-planes $\Sigma_{\pi_0,x}$ embedded in $\pi_0$ as defined
in equations \eqref{def:psub} and \eqref{def:sigmasub}
respectively.  Then, we have the following:
\begin{equation}
\text{H1}(k,P,\tilpi) \Rightarrow
\text{H1}(k-1,P_{\pi_0,x},\Sigma_{\pi_0,x})
\end{equation}
for at least half of the pairs $(\pi_0,x)$.  Furthermore, this
regularity holds at all lower levels of the induction in the sense
that (using the notation already established in this section)
\begin{equation}
\text{H1}(r,Q,\Sigma) \Rightarrow
\text{H1}(r-1,Q_{\pi_0,x},\Gamma_{\sigma_0,x})
\end{equation}
for at least half of the pairs $(\sigma_0,x)$ .
\end{lemma}

\begin{proof}
Given $(\pi_0,x) \in D(P,\tilpi)$, define sets $P_{\pi_0,x}$ and
$\Sigma_{\pi_0,x}$ as above.  Assuming that
H1($k-1$,$P_{\pi_0,x}$,$\Sigma_{\pi_0,x}$) fails when $(\pi_0,x)
\in D'(P,\tilpi)$, where $D'(P,\tilpi) \subset D(P,\tilpi)$ and
$\lvert D'(P,\tilpi) \rvert \geq \frac{1}{2}\lvert D(P,\tilpi)
\rvert$, we will obtain upper and lower bounds for the size of the
set $V_{k,p}(P,\tilpi)$ as defined in equation \eqref{def:vkp}.
For convenience, let us recall this definition:
\begin{equation*}
V_{k,p} = \{((\pi_0,\pi,p_1,\dotsc,p_k),x) \in V_k \times P: x \in P \cap (\pi \smallsetminus [p_1,\dotsc,p_k])\}
\end{equation*}
Assuming H1($k$,$P$,$\tilpi$) and H2($k$,$P$,$\tilpi$) hold, we
may use the computations from the previous section, and obtain a
lower bound of
\begin{equation}\label{h1:vkplower}
\lvert V_{k,p} \rvert \gtrapprox
\frac{\ntili^{2k+1}}{\P^{k}\ntilpi^{2k-1}}.
\end{equation}

Before proceeding, we will refine the set $V_{k,p}$ to the
following slightly smaller subset
\begin{equation}
V'_{k,p} = \{((\pi_0,\pi,p_1,\dotsc,p_k),x) \in V_{k,p}: (x,\pi_0) \in D'(P,\tilpi)\}.
\end{equation}
Since there is a uniform (throughout the set $D(P,\tilpi)$) lower
bound on the number of objects of the form $(\pi,p_1,\dotsc,p_k)
\in I_k(P,\tilpi)$ such that $((\pi_0,\pi,p_1,\dotsc,p_k),x) \in
V_{k,p}$, and $\lvert D'(P,\tilpi) \rvert \geq \frac{1}{2}\lvert
D(P,\tilpi) \rvert$, we have $\lvert V'_{k,p}(P,\tilpi)\rvert \sim
\lvert V_{k,p}(P,\tilpi) \rvert$.

Next, we compute an upper bound for the size of
$V'_{k,p}(P,\tilpi)$ in two slightly different ways.  First of
all, there are $\lesssim\P\ntilpi$ choices for the pair
$(\pi_0,x)$.  Once this pair has been fixed, we make use of the
estimate $\lvert P \cap \pi_0 \rvert \approx
\frac{\ntili}{\ntilpi}$, and observe that there are $\approx
\left(\frac{\ntili}{\ntilpi}\right)^{k-1}$ choices for the
$(k-1)$-tuple $(p_1,\dotsc,p_{k-1})$.  Next, assuming that
H1($k-1$,$P_{\pi_0,x}$,$\Sigma_{\pi_0,x}$) fails, it follows that
there are at most $O(1)$ $k$-planes from $\tilpi$ which contain
the $(k-1)$-tuple $(p_1,\dotsc,p_{k-1})$ and contribute a
$(k-1)$-plane to the set $\Sigma_{\pi_0,x}$.  This means that
there are at most $\sim \F^{k-1}$ choices for the remaining point
$p_k$.  Putting this all together, we have an estimate of
\begin{equation}\label{h1:vkpupper1}
\lvert V'_{k,p} \rvert \lessapprox
\P\ntilpi\left(\frac{\ntili}{\ntilpi}\right)^{k-1}\F^{k-1}.
\end{equation}

Also, we will need to make another, seemingly cruder estimate of
the same set.  This estimate is obtained in the same manner as
equation \eqref{h1:vkpupper1}, except we use the trivial estimate
$\lvert P \cap \pi_0 \rvert \leq \F^{k}$ in place of the estimate
$\lvert P \cap \pi_0 \rvert \approx \frac{\ntili}{\ntilpi}$.  This
yields an upper bound of
\begin{equation}\label{h1:vkpupper2}
\lvert V'_{k,p} \rvert \lessapprox \P\ntilpi\F^{k^2-1}.
\end{equation}

Combining the upper bounds \eqref{h1:vkpupper1} and
\eqref{h1:vkpupper2} with the lower bound \eqref{h1:vkplower}, we
obtain the incidence bounds
\begin{equation}
\I \lessapprox \P^{\frac{k+1}{k+2}}
\ntilpi^{\frac{k+1}{k+2}}\F^{\frac{k-1}{k+2}}
\end{equation}
and
\begin{equation}
\I \lessapprox \P^{\frac{k+1}{2k+1}}
\ntilpi^{\frac{2k}{2k+1}}\F^{\frac{k^2-1}{2k+1}}
\end{equation}
respectively.

Next, for each integer $k$, we define a quantity $\alpha (k)$ as
\begin{equation}
\alpha (k) = \frac{k^3+k^2-4k-4}{k^3+k^2-2}.
\end{equation}
While this quantity may seem a bit odd, the reader will notice
that it satisfies the inequalities $0 \leq \alpha (k) \leq 1$ for
every $k$, and (after a great deal of tedious algebra), the
following miracle occurs:
\begin{align*}
\I& \lessapprox \left(\P^{\frac{k+1}{k+2}} \ntilpi^{\frac{k+1}{k+2}}\F^{\frac{k-1}{k+2}} \right)^{\alpha (k)} \left(\P^{\frac{k+1}{2k+1}} \ntilpi^{\frac{2k}{2k+1}}\F^{\frac{k^2-1}{2k+1}} \right)^{1-\alpha (k)}\\
{}& =
\P^{\frac{k(k+1)}{k^2+2k+2}}\ntilpi^{\frac{k^2+k+2}{k^2+2k+2}}\F^{\frac{k(k+1)}{k^2+2k+2}}.
\end{align*}
As this is precisely the incidence bound that we are trying to
demonstrate, it follows that we may assume
H1($k-1$,$P_{\pi_0,x}$,$\Sigma_{\pi_0,x}$) holds for most pairs
$(\pi_0,x) \in D$.

Next, given an integer $1 \leq r \leq k-2$, we assume that
collections of points $Q$ and $(k-r+1)$-planes $\Sigma$ (as
constructed at the start of this section) satisfy hypothesis
H1($k-r+1$,$Q$,$\Sigma$):
\begin{equation}
\lvert I(Q,\Sigma) \rvert \gg \lvert Q \rvert \lvert \Sigma \rvert
^{\frac{k-r}{k-r+1}}.
\end{equation}
Suppose that H1($k-r$,$Q_{\sigma_0,x}$,$\Gamma_{\sigma_0,x}$)
fails when $(\sigma_0,x) \in D'(Q,\Sigma)$ for some set
$D'(Q,\Sigma) \subset D(Q,\Sigma)$ where$\lvert D'(Q,\Sigma)\rvert
\geq \frac{1}{2} \lvert D(Q,\Sigma)\rvert$.  then, we have the
following:
\begin{equation}\label{notbilinear:k-r}
\lvert I(Q_{\sigma_0,x},\Gamma_{\sigma_0,x})\rvert \lesssim \lvert
Q_{\sigma_0,x}\rvert \lvert \Gamma_{\sigma_0,x}
\rvert^{\frac{k-r-1}{k-r}} \text{ for $(\sigma_0,x) \in
D'(Q,\Sigma)$}.
\end{equation}

In order to show that one may assume that this situation doesn't
occur, we will proceed in a manner similar to the way we addressed
the $r=1$ case, however we will be counting objects which are a
bit more complicated that elements of the set $V_{k,p}$.

\begin{definition}
Given integers $k$ and $l$, we define a $(k,l)$-\emph{chain} to be
a $(k+2)$-tuple of points $(p_0,\dotsc,p_{k+1}) \in P^{k+2}$, and
an $l$ tuple of $k$-planes $(\pi_0,\dotsc,\pi_{l-1})\in\tilpi^l$
such that $(p_0,\dotsc,p_{k+1})$ spans a $(k+1)$-dimensional
space, and for any positive integer $m \leq l$, and any choice of
$m$ $k$-planes from the $l$-tuple $(\pi_0,\dotsc,\pi_{l-1})$, the
space $\pi_{i_1}\cap\dotsm\cap\pi_{i_m}$ is $(k-m+1)$-dimensional,
and spanned by some $(k-m+2)$-tuple of points from $(k+2)$-tuple
$(p_0,\dotsc,p_{k+1})$.  The set of all $(k,l)$-chains constructed
from the sets $P$ and $\tilpi$ will be denoted
$C_{k,l}(P,\tilpi)$.
\end{definition}

\begin{remark}
As the definition of a $(k,l)$-chain is somewhat complicated, one
might wish to regard them as $(k+1)$-simplices (with faces from
$\tilpi$ and vertices from $P$) which are missing $k-l+2$ of their
faces, but none of their vertices.
\end{remark}

This construction is motivated by the fact that when $l = r+1$,
the intersection of the $r+1$ $k$-planes in any $(k,r+1)$ chain in
the set $C_{k,r+1}(P,\tilpi)$ will be a $(k-r)$-plane from the set
$\Gamma_{\sigma_0,x}$.  Therefore, if we assume that equation
\eqref{notbilinear:k-r} holds, we will be able to demonstrate an
upper bound for $\lvert C_{k,r+1} \rvert$.  Then, after computing
a corresponding lower bound (mostly by means of Cauchy-Schwarz and
pigeonholing techniques), we will be able to show that the
incidence bound \eqref{ic:main} holds.

Our first task will be to compute a lower bound for a certain
large subclass of $C_{k,r+1}(P,\tilpi)$.  Unfortunately, this
computation is much more complicated than it was in the case
$r=1$; the reason for this is that Cauchy-Schwarz alone does not
seem to be powerful enough to count $(k,l)$-chains when $l \geq
3$.  Therefore, we shall use a ``coning procedure'', much like was
seen in the previous section for counting simplices, in order to
estimate $\lvert C_{k,r+1}(P,\tilpi)\rvert$ from below.  The basic
scheme for this construction will be to realize any
$(k,r+1)$-chain as the cone of a $(k-1,r)$-chain, which in turn is
the cone of a $(k-2,r-1)$-chain, and so on (see figure
\ref{chainconepic}).  Eventually, elements of the set
$C_{k,r+1}(P,\tilpi)$ will be described as hypercones of
$(k-r+1,2)$-chains.  Since Cauchy-Schwarz \emph{is} useful for
estimating chains of length $2$, we will arrive at the desired
lower bound for $\lvert C_{k,r+1} \rvert$, by estimating the size
of a certain class of $(k-r+1,2)$-chains, and then making a few
observations.

\begin{figure}[htb]
\begin{center}
\includegraphics{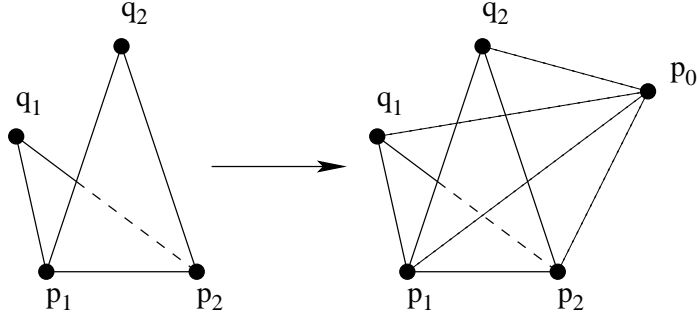}
\caption[A $(3,3)$-chain as the cone of a $(2,2)$-chain.]{This
figure illustrates the realization of a $(3,3)$-chain as the cone
of a $(2,2)$-chain.  The three $3$-planes in the $(3,3)$-chain on
the right are $\pi_0 = [p_1,p_2,q_1,q_2]$, $\pi_1 =
[p_1,p_2,p_0,q_1]$ and $\pi_2 = [p_1,p_2,p_0,q_2]$.}
\end{center}
\label{chainconepic}
\end{figure}

This computation will require a bit of extra notation.  As was
shown in the previous section, it is possible to find a collection
of disjoint pairs of points and $k$-planes $D(P,\tilpi) \subset P
\times \tilpi$ such that for each $(p_0,\pi_0) \in D(P,\tilpi)$
one has
\begin{equation}
f(p_0,\pi_0) \gtrapprox \frac{\I^{2k+1}}{\P^k\ntilpi^{2k-1}\lvert
D(P,\tilpi)\rvert},
\end{equation}
where the function $f$ is defined (for pairs $(p_0,\pi_0) \in P
\times \tilpi$) as
\begin{equation}
f(p_0,\pi_o) = \lvert \{(\pi,p_1,\dotsc,p_k)\in \tilik:
((\pi_0,\pi,p_1,\dotsc,p_k),p_0) \in V_{k,p}(P,\tilpi)\}\rvert.
\end{equation}
So, for each $(p_0,\pi_0) \in D(P,\tilpi)$, one can define
families $Q_{p_0,\pi_0}^{(k-1)}$ and $\Sigma_{p_0,\pi_0}^{(k-1)}$
of points and $(k-1)$-planes contained in $\pi_0$.  By definition,
these collections satisfy the following estimates:
\begin{align}
\lvert Q_{p_0,\pi_0}^{(k-1)} \rvert& \approx \frac{\ntili}{\ntilpi}\\
\lvert I_k(Q_{p_0,\pi_0}^{(k-1)},\Sigma_{p_0,\pi_0}^{(k-1)})&
\gtrapprox \frac{\ntili^{2k+1}}{\P^k\ntilpi^{2k-1}\lvert
D(P,\tilpi)\rvert}.
\end{align}

Next, we observe that any element of the set $C_{k,r+1}(P,\tilpi)$
which contains the point $p_0$ as a vertex and the $k$-plane
$\pi_0$ as a face must be the cone of a unique element of the set
$C_{k-1,r}(Q_{p_0,\pi_0}^{(k-1)},\Sigma_{p_0,\pi_0}^{(k-1)})$.
This observation allows us to conclude that
\begin{equation*}
\lvert C_{k,r+1}(P,\tilpi)\rvert \gtrapprox \sum_{(p_0,\pi_0)\in
D(P,\tilpi)}\lvert
C_{k-1,r}(Q_{p_0,\pi_0}^{(k-1)},\Sigma_{p_0,\pi_0}^{(k-1)})\rvert.
\end{equation*}
So, by a simple application of the pigeonhole principle, we may
choose a single value $(p_0,\pi_0) \in D(P,\tilpi)$ such that
\begin{equation}
\lvert C_{k,r+1}(P,\tilpi)\rvert \gtrapprox \lvert
C_{k-1,r}(Q_{p_0,\pi_0}^{(k-1)},\Sigma_{p_0,\pi_0}^{(k-1)})\rvert
\lvert D(P,\tilpi)\rvert.
\end{equation}

Therefore, we have reduced our problem to that of estimating the
size of the sets
$C_{k-1,r}(Q_{p_0,\pi_0}^{(k-1)},\Sigma_{p_0,\pi_0}^{(k-1)})$.  At
this point, the procedure repeats itself, this time with the
collections $Q_{p_0,\pi_0}^{(k-1)}$ and
$\Sigma_{p_0,\pi_0}^{(k-1)}$ of points and $(k-1)$-planes.

\begin{remark}
Since these estimates are hold uniformly for a specifically chosen
(via the pigeonhole principle) pair $(p_0,\pi_0) \in D(P,\tilpi)$,
we will no longer use these parameters as subscripts when
describing induced sets of points or $(k-1)$-planes.  For example,
the set $\Sigma^{(k-1)}$ will be used to denote an induced family
of $(k-1)$-planes relative to some generic $(p_0,\pi_0) \in
D(P,\tilpi)$.  This convention is purely in the interest of
simplifying notation, and will be used for analogous constructions
to appear later in the proof.  Furthermore, the sets
$\Gamma_{\sigma_0,x}$ referred to in equation
\eqref{notbilinear:k-r} will henceforth be denoted
$\Sigma^{(k-r)}$.
\end{remark}

In general, once the sets $Q^{(k-l)}$ and $\Sigma^{(k-l)}$ have
been constructed (for some integer $1 \leq l \leq r-2$), one can
choose a collection $D(Q^{(k-l)},\Sigma^{(k-l)}) \subset
Q^{(k-l)}\times\Sigma^{(k-l)}$ of disjoint pairs of points and
$(k-l)$-planes, and use this collection to define (as was
described at the beginning of this section) to create families
$Q^{(k-l-1)}$ and $\Sigma^{(k-l-1)}$ of points and
$(k-l-1)$-planes which satisfy the estimates
\begin{align}
\lvert Q^{(k-l-1)}\rvert& \approx \frac{\lvert I(Q^{(k-l)},\Sigma^{(k-1)})\rvert}{\lvert\Sigma^{(k-l)}\rvert}\label{j:inductive1}\\
\lvert I_{k-l}(Q^{(k-l-1)},\Sigma^{(k-l-1)})\rvert& \gtrapprox
\frac{\lvert I(Q^{(k-l)},\Sigma^{(k-l)})\rvert^{2(k-l)+1}}{\lvert
Q^{(k-l)} \rvert^{k-l}\lvert \Sigma^{(k-l)}
\rvert^{2(k-l)-1}\lvert
D(Q^{(k-l)},\Sigma^{(k-l)})\rvert}\label{j:inductive2}.
\end{align}
Furthermore, by identifying elements of the sets
$C_{k-l,r+1-l}(Q^{(k-l)},\Sigma^{(k-l)})$ and
$C_{k-l-1,r-l}(Q^{(k-l-1)},\Sigma^{(k-l-1)})$, we have the
estimate
\begin{equation}
\lvert C_{k-l,r+1-l}(Q^{(k-l)},\Sigma^{(k-l)}) \rvert \gtrapprox
\lvert C_{k-l-1,r-l}(Q^{(k-l-1)},\Sigma^{(k-l-1)}) \rvert \lvert
D(Q^{(k-l)},\Sigma^{(k-l)}) \rvert.
\end{equation}
Since these observations hold for any integer $1 \leq l \leq r-2$,
we arrive at the following conclusion:
\begin{equation}\label{crk:prod}
\lvert C_{k,r+1}(P,\tilpi)\rvert \gtrapprox \lvert
C_{k-r+1,2}(Q^{(k-r+1)},\Sigma^{(k-r+1)})\rvert \lvert D(P,\tilpi)
\rvert \prod_{j=k-r+2}^{k-1}\lvert D(Q^{(j)}\Sigma^{(j)}) \rvert.
\end{equation}

At this point, since we have reduced to problem to estimating the
size of a family of $(k-r+1,2)$-chains, we may use Cauchy-Schwarz.
The construction of elements of the set
$C_{k+r-1,2}(Q^{(k-r+1)},\Sigma^{(k-r+1)})$ is virtually identical
to the construction of the sets
$V_{2,p}(Q^{(k-r+1)},\Sigma^{(k-r+1)})$, except that we add an
extra point to each of the $(k-r+1)$-planes, rather than just one
of them.  Therefore, Cauchy-Schwarz yields an estimate of
\begin{align}
\lvert C_{k+r-1,2}(Q^{(k-r+1)},\Sigma^{(k-r+1)}) \rvert & \gtrsim \frac{\lvert I_{k-r+1}(Q^{(k-r+1)},\Sigma^{(k-r+1)})\rvert^2}{\lvert Q^{(k-r+1)}\rvert^{k-r+1}}\notag\\
{}& \quad \frac{\lvert I(Q^{(k-r+1)},\Sigma^{(k-r+1)})\rvert^2}{\lvert \Sigma^{(k-r+1)}\rvert^2}\notag\\
{}& \gtrsim \frac{\lvert
I(Q^{(k-r+1)},\Sigma^{(k-r+1)})\rvert^{2(k-r+2)}}{\lvert
Q^{(k-r+1)}\rvert^{k-r+1} \lvert
\Sigma^{(k-r+1)}\rvert^{2(k-r+1)}}.
\end{align}

Furthermore, since we are assuming (by induction) that the
hypothesis H1($k-r+1$,$Q^{(k-r+1)}$,$\Sigma^{(k-r+1)}$) holds, it
follows that these $(k-r+1,2)$-chains are non-degenerate (i.e. the
two $(k-r+1)$-planes contained in each are distinct).  Explicitly,
the degenerate $(k-r+1)$-chains which we would like to discard are
simply elements of the set
$I_{k-r+1}(Q^{(k-r+1)},\Sigma^{(k-r+1)})$ equipped with an
additional pair of points, so it has size $\approx \lvert
I_{k-r+1}(Q^{(k-r+1)},\Sigma^{(k-r+1)})\rvert \frac{\lvert
I_{k-r+1}(Q^{(k-r+1)},\Sigma^{(k-r+1)})\rvert^2}{\lvert
\Sigma^{(k-r+1)}\rvert^2}$.  Applying the hypothesis
H1($k-r+1$,$Q^{(k-r+1)}$,$\Sigma^{(k-r+1)}$) yields
\begin{align*}
\lvert C_{k+r-1,2}(Q^{(k-r+1)},\Sigma^{(k-r+1)}) \rvert &\gtrsim \frac{\lvert I_{k-r+1}(Q^{(k-r+1)},\Sigma^{(k-r+1)})\rvert^2}{\lvert Q^{(k-r+1)}\rvert^{k-r+1}}\\
{}& \quad \frac{\lvert I(Q^{(k-r+1)},\Sigma^{(k-r+1)})\rvert^2}{\lvert \Sigma^{(k-r+1)}\rvert^2}\\
{}& \gtrsim \frac{\lvert I_{k-r+1}(Q^{(k-r+1)},\Sigma^{(k-r+1)})\rvert^{k-r+1}}{\lvert Q^{(k-r+1)}\rvert^{k-r+1}\lvert \Sigma^{(k-r+1)}\rvert^{k-r}}\\
{}& \quad\lvert I_{k-r+1}(Q^{(k-r+1)},\Sigma^{(k-r+1)})\rvert\\
{}& \quad \frac{\lvert I(Q^{(k-r+1)},\Sigma^{(k-r+1)})\rvert^2}{\lvert \Sigma^{(k-r+1)}\rvert^2}\\
{}& \gg \lvert I_{k-r+1}(Q^{(k-r+1)},\Sigma^{(k-r+1)})\rvert \\
{}& \quad\frac{\lvert
I(Q^{(k-r+1)},\Sigma^{(k-r+1)})\rvert^2}{\lvert
\Sigma^{(k-r+1)}\rvert^2}.
\end{align*}
Therefore, we have
\begin{multline}
\lvert C_{k,r+1}(P,\tilpi)\rvert \gtrapprox\\
\frac{\lvert
I(Q^{(k-r+1)},\Sigma^{(k-r+1)})\rvert^{2(k-r+2)}}{\lvert
Q^{(k-r+1)}\rvert^{k-r+1} \lvert
\Sigma^{(k-r+1)}\rvert^{2(k-r+1)}}\lvert D(P,\tilpi) \rvert
\prod_{j=k-r+2}^{k-1}\lvert D(Q^{(j)}\Sigma^{(j)}) \rvert.
\end{multline}

Before preceding, we must refine the set
$C_{k-r+1,2}(Q^{(k-r+1)},\Sigma^{(k-r+1)})$ to a slightly smaller
set $C'_{k-r+1,2}(Q^{(k-r+1)},\Sigma^{(k-r+1)})$ defined as
\begin{multline}
C'_{k-r+1,2}(Q^{(k-r+1)},\Sigma^{(k-r+1)}):=\\
 \left\{
\begin{array}{l}
(\sigma_0,\sigma,q_0,\dotsc,q_{k-r+2})\\
\in C_{k-r+1,2}(Q^{(k-r+1)},\Sigma^{(k-r+1)})
\end{array}:
\begin{array}{l}
(\sigma_0,q_0)\\
\in D'(Q^{(k-r+1)},\Sigma^{(k-r+1)})
\end{array}
\right\}.
\end{multline}
Making use of the uniform lower bound estimate from equation
\eqref{frestimate}, and the fact that $\lvert
D'((Q^{(k-r+1)},\Sigma^{(k-r+1)})\rvert \geq \frac{1}{2}\lvert
D((Q^{(k-r+1)},\Sigma^{(k-r+1)})\rvert$, we have
\begin{equation*}
\lvert C'_{k-r+1,2}(Q^{(k-r+1)},\Sigma^{(k-r+1)})\rvert \sim
\lvert C_{k-r+1,2}(Q^{(k-r+1)},\Sigma^{(k-r+1)})\rvert.
\end{equation*}
We shall denote the elements of the set $C_{k,r+1}(P,\tilpi)$
which are cones of $(k-r+1,2)$-chains from
$C_{k-r+1,2}(Q^{(k-r+1)},\Sigma^{(k-r+1)})$ as
$C'_{k,r+1}(P,\tilpi)$.

Next, it is necessary to simplify the expression on the right hand
side this expression.  Making use of equations
\eqref{j:inductive1} and \eqref{j:inductive2}, we have
\begin{multline}
\frac{\lvert I(Q^{(j_0-1)},\Sigma^{(j_0-1)})\rvert^{j_0(j_0-k+r)}}{\lvert Q^{(j_0-1)}\rvert^{(j_0-1)(j_0-k+r-1)+j_0-k+r-2}\lvert \Sigma^{(j_0-1)}\rvert^{(j_0-1)(j_0-k+r)}}\\
 \lvert D(P,\tilpi)\rvert \prod_{j=j_0}^{k-1}\lvert D(Q^{(j)},\Sigma^{(j)})\lvert \\
\gtrapprox \left(\frac{\lvert \Sigma^{(j_0)}\rvert}{\lvert I(Q^{(j_0)},\Sigma^{(j_0)})\rvert}\right)^{j_0(j_0-k+r-1)-1}\left(\frac{\lvert I(Q^{(j_0)},\Sigma^{(j_0)})\rvert^{2j_0+1}}{\lvert Q^{(j_0-1)} \rvert^{j_0}\lvert \Sigma^{(j_0)}\rvert^{2j_0-1}\lvert D(Q^{(j_0)},\Sigma^{(j_0)}) \rvert}\right)^{j_0-k+r}\\
\lvert D(P,\tilpi)\rvert \prod_{j=j_0}^{k-1}\lvert
D(Q^{(j)},\Sigma^{(j)})\lvert
\end{multline}
holds when $k-r+2 \leq j_0 \leq k-1$.  Furthermore, making use of
the trivial estimate
\begin{equation*}
\lvert D(Q^{(j_0)},\Sigma^{(j_0)}) \rvert \leq \lvert Q^{(j_0)}
\rvert \lvert \Sigma^{(j_0)} \rvert,
\end{equation*}
we arrive at an estimate of
\begin{multline}
\lvert C'_{k,r+1}(P,\tilpi)\rvert \gtrapprox\\
\frac{\lvert I(Q^{(j_0)},\Sigma^{(j_0)})\rvert^{(j_0+1)(j_0-k+r+1)}}{\lvert Q^{(j_0)}\rvert^{j_0(j_0-k+r-1)+j_0-k+r-1}\lvert \Sigma^{(j_0)}\rvert^{j_0(j_0-k+r+1)}}\\
 \lvert D(P,\tilpi)\rvert \prod_{j=j_0+1}^{k-1}\lvert D(Q^{(j)},\Sigma^{(j)})\lvert.
\end{multline}

After repeating this procedure $r-2$ times, and making use of the
trivial estimate $\lvert D(P,\tilpi) \rvert \leq \P \ntilpi$, we
finally arrive at the desired conclusion
\begin{equation}
\lvert C'_{k,r+1}(P,\tilpi)\rvert \gtrapprox \frac{\lvert I(P,\tilpi)\rvert^{(k+1)(r+1)}}{\P^{kr+k-1}\ntilpi^{k(r+1)}}\\
\end{equation}

Next, we shall compute a complementary upper bound for $\lvert
C'_{k,r+1}(P,\tilpi) \rvert$ in two different ways (as was done
for the case $r=1$).  First, choose a pair $(p_0,\pi_0) \in P
\times \tilpi$ such that $p_0 \notin \pi_0$; there are clearly at
most $\P\ntilpi$ choices for this pair.  By definition, the
remaining $k+1$ points must lie inside $\pi_0$, and they will
uniquely define the remaining $r$ $k$-planes once chosen.  So,
once the number of ways to choose the remaining points has been
computed, we will have an upper bound for $\lvert
C'_{k,r+1}(P,\tilpi) \rvert$.  Since $\lvert P \cap \pi_o \rvert
\approx \frac{\I}{\ntilpi}$, and these points cannot cluster on
low dimensional subspaces of $\pi_0$, there are $\approx
\left(\frac{\I}{\ntilpi}\right)^{k-r}$ choices for the first $k-r$
points of the $(k,r+1)$-chain.  Next, assuming equation
\eqref{notbilinear:k-r} holds, there are at most $\sim 1$ choices
for the $(k-r)$-plane found at the intersection of the
$(r+1)$-fold intersection of the $k$-planes in the $(k,r+1)$
chain.  This means that we have at most $\sim \F^{k-r}$ choices
for the next point.  Finally, there are $\approx
\left(\frac{\I}{\ntilpi}\right)^r$ choices for the remaining $r$
points.  These observations yield the bound
\begin{equation}\label{ckrupper1}
\lvert C'_{k,r+1}(P,\tilpi) \rvert \lessapprox
\frac{\P\I^k\F^{k-r}}{\ntilpi^{k-1}}.
\end{equation}

Next, we make the same computation, only use the cruder estimate
$\lvert P \cap \pi_0 \rvert \leq F^k$ instead of $\lvert P \cap
\pi_0 \rvert \approx \frac{\I}{\ntilpi}$.  This method yields the
bound
\begin{equation}\label{ckrupper2}
\lvert C'_{k,r+1}(P,\tilpi) \rvert \lessapprox \P \ntilpi
\F^{k^2+k-r}.
\end{equation}

In order to finish the proof, we simply combine the upper and
lower bounds we have now computed for $\lvert C'_{k,r+1}(P,\tilpi)
\rvert$.  Making use of the result from the previous computation,
and the estimates \eqref{ckrupper1} and \eqref{ckrupper2}, we have
\begin{equation}
\lvert I(P,\tilpi) \rvert \lessapprox
\P^{\frac{r(k+1)}{(k+1)(r+1)-k}}\ntilpi^{\frac{(k+1)(r+1)-k-r}{(k+1)(r+1)-k}}\F^{\frac{k-r}{(k+1)(r+1)-k}},
\end{equation}
and
\begin{equation}
\lvert I(P,\tilpi) \rvert \lessapprox
\P^{\frac{r(k+1)}{(k+1)(r+1)}}\ntilpi^{\frac{(k+1)(r+1)-r}{(k+1)(r+1)-k}}\F^{\frac{k^2+k-r}{(k+1)(r+1)}}
\end{equation}
respectively.

As it turns out, our main estimate \eqref{ic:main} is simply a
convex combination of these two estimates.  To be more specific,
if we define the quantity $\alpha_r(k)$ as
\begin{equation}
\alpha_r(k) =
\left(\frac{k(k+1)(r+1)-r(k^2+2k+2)}{k^2+2k+2}\right)\left(\frac{(k+1)(r+1)-k}{kr}\right),
\end{equation}
then, one can compute that $0 \leq \alpha_{r}(k) \leq 1$ for
appropriate values of $k$ and $r$, and
\begin{align*}
\lvert I(P,\tilpi) \rvert& \lessapprox (\P^{\frac{r(k+1)}{(k+1)(r+1)-k}}\ntilpi^{\frac{(k+1)(r+1)-k-r}{(k+1)(r+1)-k}}\F^{\frac{k-r}{(k+1)(r+1)-k}})^{\alpha_r(k)}\\
{}& \quad (\P^{\frac{r(k+1)}{(k+1)(r+1)}}\ntilpi^{\frac{(k+1)(r+1)-r}{(k+1)(r+1)-k}}\F^{\frac{k^2+k-r}{(k+1)(r+1)}})^{1-\alpha_r(k)}\\
{}& =
\P^{\frac{k(k+1)}{k^2+2k+2}}\ntilpi^{\frac{k^2+k+2}{k^2+2k+2}}\F^{\frac{k(k+1)}{k^2+2k+2}}.
\end{align*}
\end{proof}


\section{Simplex Construction Part Three: The Upper Bound}\label{upperbound}

Now that we have established a lower bound on the number of
simplices appearing in the configuration of points and $k$-planes,
we must find an upper bound. This task requires making a few
observations about the set $V_k$.  First of all, by construction,
the pair of $k$-planes in any element of $v \in V_k$ spans some
affine $(k+1)$-dimensional space $\Lambda_v$.  Since our original
family of $k$-planes (and hence every refinement of that family)
is direction separated, we have
\begin{equation}
\lvert \{\pi \in \tilpi:\pi \subset \Lambda_v\} \rvert \lesssim
\F^k
\end{equation}
for all $v \in V_k$.  So, when constructing a simplex from some $v
\in V_k$ by choosing the remaining $k$ faces from $\Lambda_v$, we
have $\F^k$ choices for each face.

Before going ahead with this construction, however, one must
observe that an element $v \in V_k$ is not merely a pair of
$k$-planes; each such $v$ comes equipped with a $k$-tuple of
points on its ``spine'' (the intersection of the two $k$-planes).
In order to manage this technicality, recall that (from the
refined definition of $\tilik$), the spine of any $v \in V_k$ must
carry at least $\frac{\ntili}{10\ntilpi\F}$ points from $P$.  So,
if we delete the points from each element of $V_k$ to create the
following set
\begin{equation}
V_{k,del}:= \left\{(\pi_0,\pi)\in \tilpi \times \tilpi:
\begin{array}{l}
\exists (p_1,\dotsc,p_k) \in P^k \text{ so that}\\
(\pi_0,\pi,p_1\dotsc,p_k) \in V_k
\end{array} \right\},
\end{equation}
we have the following bound:
\begin{equation}
\lvert V_k \rvert \gtrsim \left(\frac{\ntili}{\ntilpi\F}\right)^k
\lvert V_{k,del}\rvert.
\end{equation}

Now that we have dealt with this possible over-counting, we may
construct simplices from each element $v \in V_{del,k}$ by
choosing the $k$ remaining faces from $\Lambda_v$.  This yields an
upper bound in the number of simplices of
\begin{equation}
\lvert S_k \rvert \lessapprox \nvk
\left(\frac{\ntilpi\F}{\ntili}\right)^k\F^{k^2}.
\end{equation}

If we combine this with the lower bound obtained in Proposition
\ref{lowerlemma}, we arrive at the desired incidence bound:
\begin{equation}
\ntili \lessapprox
\P^{\frac{k(k+1)}{k^2+2k+2}}\ntilpi^{\frac{k^2+k+2}{k^2+2k+2}}\F^{\frac{k(k+1)}{k^2+2k+2}}.
\end{equation}


\end{document}